\documentclass{amsart}
\usepackage{amsmath,amsthm,amssymb,amsfonts}
\usepackage[colorlinks=true]{hyperref}
\usepackage{eucal}

\hypersetup{urlcolor=blue, citecolor=blue}

\let\MR\mr

\def\doi#1{
   {\href{http://dx.doi.org/#1}
   {{\mdseries\ttfamily DOI}}}}

\newcommand{\al}{\alpha}    \newcommand{\be}{\beta}
\newcommand{\de}{\delta}    
\newcommand{\vep}{\epsilon}  \newcommand{\ep}{\varepsilon}

\newcommand{\ga}{\gamma}    \newcommand{\Ga}{\Gamma}
\def\<{\langle}             \def\>{\rangle}
\newcommand{\R}{\mathbb{R}}
\newcommand{\N}{\mathbb{N}}

\newcommand{\pt}{\partial_t}\newcommand{\pa}{\partial}

\newcommand{\beeq}{\begin{equation}}\newcommand{\eneq}{\end{equation}}

\newenvironment{prf}{\noindent {\bf Proof.} }{\endprf\par}
\def \endprf{\hfill  {\vrule height6pt width6pt depth0pt}\medskip}

\numberwithin{equation}{section}

\newcommand{\gm}{\mathfrak{g}}

\newcommand\od{\kappa}
\newcommand{\eps}{\varepsilon}

\def\I{{\mathcal{I}}}

\def\X{{\mathcal{X}}}

\def\<{\langle}
\def\>{\rangle}
\def\({\left(}
\def\){\right)}

\newtheorem{lemma}{Lemma}[section]
\newtheorem{theorem}{Theorem}[section]

\theoremstyle{remark}
\newtheorem{remark}{Remark}

\theoremstyle{definition}

\title[Quasilinear Wave Equations]
{Global existence of null-form wave equations on small asymptotically Euclidean manifolds}

\author{Chengbo Wang}
\address{Department of Mathematics\\
                Zhejiang University\\
                Hangzhou 310027, China}
\email{wangcbo@gmail.com}
\urladdr{http://www.math.zju.edu.cn/wang/}
\thanks{The first author was supported by Zhejiang Provincial Natural Science Foundation of China LR12A01002, the Fundamental Research Funds for the Central Universities, NSFC 11301478, 11271322 and J1210038.}

\author{Xin Yu}
\email{yumei165@gmail.com}

\date{}
\dedicatory{} \commby{}

\begin{document}

\begin{abstract}{We prove the global existence of the small solutions to the Cauchy problem for quasilinear wave equations satisfying the null condition on $(\R^3,\gm)$,
where the metric $\gm$ is a small perturbation of the flat metric
and approaches the Euclidean metric like $(1+|x|^2)^{-\rho/2}$ with $\rho>1$. Global and almost global existence for systems without the null condition are also discussed for certain small time-dependent perturbations of the flat metric in the appendix.
}
\end{abstract}

\keywords{
null condition; asymptotically Euclidean manifold; asymptotically flat manifold; KSS estimates; local energy estimates; almost global existence}

\subjclass[2010]{35L72, 58J45, 35L05}

\maketitle 
\numberwithin{equation}{section}

\section{Introduction}

This paper is concerned with the Cauchy problem for the system of the quasi-linear wave equations, with multiple speeds, in three space
dimensions of the form
\beeq
\begin{cases}
\partial_t^2 u^I-c^2_I \triangle_\gm u^I
=Q^{I,\al\be\ga}_{JK}\partial_\al u^J
\partial_\be\partial_\ga u^{K}+S^{I,\al\be}_{JK}\partial_\al u^J
\partial_\be u^{K}
\\
u^I(0,x)=u^I_0, \pa_t u^I(0,x)=u^I_1, I=1,\cdots,M
\end{cases}
\eneq
subject to suitably small initial conditions, posed on certain asymptotically Euclidean manifolds $(\R^3,\gm)$.
Here, we will assume that the speeds of propagation $c_I$ are positive and distinct, and we refer this situation
as the nonrelativistic case. Straightforward modifications of the argument will give the more general case where the various components are allowed to have the same speed. It will be apparent that our argument apply to the general system with quadratic (with null condition) and higher order (cubic and higher) terms. To simplify the presentation, we restrict ourselves to the case of quadratic level of perturbation.

We shall construct a unique global classical solution, provided that the coefficients of the nonlinear terms
satisfy the null condition and the metric $\gm$ is a small perturbation of the flat Euclidean metric
and approaches the Euclidean metric like $(1+|x|^2)^{-\rho/2}$ with $\rho>1$.

In the case of flat metric $\gm$, this problem has been
extensively studied. When all of the speeds are the same, the null condition was first identified by Klainerman and
shown to have global existence of small solutions in \cite{Ch}
and \cite{Kl2} (see also \cite{Hormander, J2}).
Without the null condition, we can only have almost global existence in general (see \cite{JK, Kl1, J2} for existence and \cite{J, Sid1} for blow up results). Notice that this does not mean that the null condition is the necessary condition for the general quadratic quasilinear problems to admit global solutions with small data. Actually, there is a larger class of nonlinearities which will ensure global existence, which is related to the so-called ``weak null condition" of Lindblad and Rodnianski, see e.g. \cite{LdRo03, LdRo05, Al06, Ld08, KaKu08-3}.

Small solutions always have global existence in higher dimensions \cite{KlPon, Sh, Kl1}.
The two-dimensional case is rather more complicated. The sharp results are given in \cite{Al2, Al1}, with previous works in \cite{Go93, Ho, Ka}.

In the nonrelativistic case, the null condition guarantees that the self-interaction of each wave family is nonresonant
and is the natural one for systems of quasilinear wave equations with multiple
speeds. It is equivalent to the requirement that no plane wave solution of the system is
genuinely nonlinear (see e.g. \cite{Jo90, AgYo98}).
By now there is an extensive literature devoted to this topic; without being exhaustive we
mention \cite{Kov, AgYo98, Kubo, yo, ST, KuYo01, So03, So08, Hi04, Ka04, Ka05, KaYo06, KaKu08-2}. It is remarkable that the approach of \cite{ST} and \cite{Hi04} does not use estimates of the fundamental solution for the free wave equation, which seems more robust when considering problems with variable coefficients or problems in exterior domains, see e. g. \cite{MeSo07}. The two-dimensional case has also been examined, see \cite{Ho1, Kubo, Hi04} and references therein.

In exterior domains, null form quasilinear wave equations were previously
studied in \cite{KSS02-2, MeSo05, MeSo07, MeNaSo05, MeNaSo05-2, KaKu08}. The general quasilinear problems
(without null conditions) 
 were also studied, see \cite{KSS04, MeSo06} and references therein.

It is interesting to investigate similar problems on various space-time manifolds. Recently, there have been some progress in this direction. Global and almost global existence of the solutions for the 
semilinear problems posed on asymptotically Euclidean manifolds have been obtained in \cite{BoHa}, \cite{SoWa10} and \cite{WaYu11}. Global existence of the solutions for the null form semilinear wave equations on slowly rotating Kerr spacetimes or time dependent inhomogeneous media (compact perturbation) has been obtained in \cite{Lu} and \cite{Ya}. In this paper, we will deal with the nonrelativistic quasilinear wave equations on small asymptotically Euclidean manifolds, mainly inspired by the approach of \cite{ST, Hi04}, together with the local energy estimates with variable coefficients obtained in \cite{HWY10}
(see also \cite{MeSo06} and references therein). To our knowledge, our result is the first work studying the global existence for quasilinear wave equations on asymptotically flat manifolds. It will be interesting if we can deal with the problems with general non-trapping asymptotically Euclidean manifolds or small asymptotically flat manifolds with certain time-dependent metrics.

Before stating our main result, we introduce the necessary notations.
Points in $\R^4$ will be denoted by $(x^0,x^1,x^2,x^3)=(t,x)$.
Partial derivatives will be written as
$ \partial_\al={\partial}/{\partial x^\al}$,
with the abbreviations
$ \partial=(\partial_0,\partial_1,\partial_2,\partial_3)=(\partial_t,\nabla)$.
Here, we have used the convention that Greek indices range from $0$ to $3$ and Latin indices from $1$ to $3$. We will also abuse the notation to use the Greek indices to denote multi-indices, which should be clear in the context.
Hereafter, the Einstein summation convention will be performed over repeated indices.
The rotational vector fields are defined as
\begin{equation*}
\Omega_{ij}=x^i \pa_j-x^j \pa_i, 1\le i<j\le 3,
\end{equation*} and
the scaling vector field is defined by
\begin{equation}
\label{scaling}
S=t\partial_t+r\partial_r=x^\alpha\partial_\alpha,\quad r=|x|=\sqrt{(x^1)^2+(x^2)^2+(x^3)^2}.
\end{equation}
The collection of these eight vector fields will be labeled as
\[
\Gamma=(\Gamma_0,\ldots,\Gamma_7)=(\partial,\Omega,S).
\]

 We consider the asymptotically
Euclidean Riemannian manifolds $( \R^3 , \gm)$ with
\begin{equation*}
\gm =  g_{ij} (x) \, d x^i \, d x^j .
\end{equation*}
The metric $\gm$ is assumed to be a small perturbation of the flat metric. More precisely, we suppose $g_{ij} (x) \in C^{\infty} ( \R^{3} )$ and, for some
fixed $\rho
>1$ and $\delta\ll 1$,
\begin{equation}\tag{H1} \label{H1}
\forall \al \in \N^3 \qquad |\partial^\al_x ( g_{ij} -
\delta_{ij} )|\le C_\al \delta  \< x \>^{- \vert \al \vert - \rho}  ,
\end{equation}
with $\delta_{ij}=\delta^{ij}$ being the Kronecker delta function and $\<x\>=\sqrt{1+|x|^2}$. Since $\delta\ll 1$, it is clear that the metric $\gm$ is a non-trapping perturbation. Let $g = \det (g_{ij})$, the Laplace--Beltrami
operator associated with $\gm$ is given by
$$\Delta_{\gm} =  \sqrt{g}^{-1} \partial_i g^{ij} \sqrt{g} \partial_j ,
$$
where $(g^{ij} (x))$ denotes the inverse matrix of $(g_{ij}(x))$.

Consider the initial value problem for the nonlinear equations of the form
\begin{equation}
\label{pde}
(\Box_\gm u)^I \equiv (\pt^2-c_I^2 \Delta_\gm) u^I = N^I(u,u), I=1,2,\cdots, M
\end{equation}
in which the  quadratic nonlinearity $N=Q+S$ is of the form
\begin{equation}
\label{non}
Q^I(u,v)= Q^{I,\al\be\ga}_{JK}\partial_\al u^J
\partial_\be\partial_\ga v^{K},\, S^I(u,v)=S^{I,\al\be}_{JK}\partial_\al u^J
\partial_\be v^{K}.
\end{equation}
The construction of solutions will depend on the energy integral method, which requires the quasilinear part
to be symmetric:
\begin{equation}
\label{symmetric}
Q^{I,\al\be\ga}_{JK}=Q^{I,\al\ga\be}_{JK}=Q^{K,\al\be\ga}_{JI}.
\end{equation}
The key assumption  for global existence is the following
{\em null} condition which says that the self-interaction
of each wave family is nonresonant:
\begin{equation}
\label{nc}
Q^{I,\al\be\ga}_{II}\xi_\al \xi_\be \xi_\ga=
S^{I,\al\be}_{II}\xi_\al \xi_\be =0
\quad\mbox{for all}\ \xi\ \mbox{s.t.}\quad \xi_0^2=c_I^2(\xi_1^2+\xi_2^2+\xi_3^2)\ .
\end{equation}

The standard energy norm is denoted as
\[E_1(u(t))=\frac 12 \sum_{I=1}^M\int_{\R^3}
 |\pa u^I (t,x)|^2  dx,\]
and higher order derivatives will be estimated through
\begin{equation}
\label{ennorm}
E_m(u(t))=\sum_{|\al|\le m-1} E_1(\Gamma^\al u(t)),
\qquad m=2, 3, \cdots\ .
\end{equation}

In order to describe the solution space, we introduce the
time-independent vector fields
$\Lambda=(\Lambda_1,\ldots,\Lambda_7)=(\nabla,\Omega,r\partial_r)$.
Define
\[
H^m_\Lambda(\R^3)=\{f\in L^2(\R^3) : \Lambda^\al f\in L^2,
\; |\al|\le m\},
\]
with the norm
\begin{equation}
\label{sobnorm}
\|f\|_{H^m_\Lambda}=\sum_{|\al|\le m}\|\Lambda^\al f\|_{L^2}.
\end{equation}

Solutions will be constructed in the space
$\dot{H}_\Gamma^m(T)$, which is the closure of $C^\infty([0,T);C_0^\infty(\R^3))$
with respect to the norm $\sup\limits_{0\le t< T} E_m^{1/2}(u(t))$.
Thus,
$$
\dot{H}_\Gamma^ m(T)\subset\left\{u(t,x) : \partial u(t,\cdot)
\in\bigcap_{j=0}^{ m-1} C^j([0,T);H^{ m-1-j}_\Lambda)\right\}.
$$
By Sobolev embedding, it follows that
$\dot{H}_\Gamma^m(T)\subset C^{m-2}([0,T)\times\R^3)$.

An important intermediate role will be played by the following two weighted norms
\begin{equation}
\label{wnorm}
\X_ m(u(t))= \sum_{I=1}^M\sum_{|\al|=2}\sum_{|\be|\le m-2}
\|\langle c_I t-r\rangle\partial^\al\Gamma^\be u^I (t)\|_{L^2(\R^3)},
\end{equation}
and
\begin{equation}
\label{inorm}
\I_m(u(t))= \sum_{I=1}^M\sum_{|\al|\le m-1}
\left\|r^{-1/2+\mu}\<r\>^{-\mu'} \left(|\partial \Gamma^\al u^I (t)|+\frac{\Gamma^\al u^I (t)}r\right)\right\|_{L^2(\R^3)}^2
\end{equation}
with $\mu\in (0,1/2)$ and $\mu'>\mu$ to be determined later (the choice will be $\mu=1/4$ and $\mu'=\min(2\rho-1,3)/4$, see \eqref{120630-1}, \eqref{120630-2}).
The second norm is extracted from the local energy norm (also known as KSS-type estimate), which is defined as
\begin{equation}
\label{LEnorm}
LE_m(t)= \int_0^t \I_m(u(\tau))d\tau\ .
\end{equation}
 In the case of $\mu'=\mu$, such norm will be the KSS norm and we will denote such norm by $KSS_m(t)$.

Let us now state our main result precisely.
\begin{theorem}
\label{thm}
Let $\rho>1$ and $\de\ll 1$.
Assume that the nonlinear terms in \eqref{non}
satisfy the symmetric and null conditions \eqref{symmetric}, \eqref{nc}.
Then there exist constants $\eps\ll 1\ll C_0$, such that the Cauchy problem for \eqref{pde} has a unique global solution $u\in \dot H_\Gamma^\od(T)$ for
every $T>0$, when the initial data
\[
\partial u(0)\in H^{\od-1}_\Lambda(\R^3),\quad \od\ge 9
\]
satisfying
\begin{equation}
\label{smallness}
E_{\od-2}^{1/2}(u(0))\;\exp\;C_0 E_\od^{1/2}(u(0)) < \eps.
\end{equation}
Moreover, the solution satisfies the bounds for some $C_1$,
$C_2\ge1$,
\[
LE_{\od-2}(t)+E_{\od-2}(u(t))<2 C_1 \eps\quad\mbox{and}
\quad LE_{\od}(t)+E_\od(u(t))\le C_2 E_\od(u(0)) \<t\>^{C_2{\eps}}.
\]
\end{theorem}

To conclude the introduction, let us give some remarks and comments.
\begin{remark}
We remark here that, the situation for the null form quasilinear problems seems much more delicate than the general quasilinear problem, technically due to the occurrence of the scaling vector field in our argument. Actually, it is not hard to see that we can prove the almost global existence (and global existence for higher dimension) for the solutions to the quasilinear quadratic equations, on asymptotically flat manifolds with small time-dependent metric perturbation, by combining our argument with the approach in \cite{MeSo06}, where there is no need to use the scaling vector field. See Appendix \ref{sec-App} for the proof. Although there is work \cite{KaKu08-2} dealing with null form problems without using the scaling vector field in the literature, it seems difficult to be adapted for the setting of time-dependent asymptotically flat perturbation (except the case of compact perturbation, see \cite{KaKu08}).
\end{remark}

\begin{remark}
The case $n=2$ seems more difficult to handle, partly because of the lack of the local energy estimates with variable coefficients.
\end{remark}

\begin{remark} The same argument can yield global results for the system with repeated speeds, by strengthening the null condition to be nonresonant
interaction among the  waves with  the same wave speeds. See \cite{ST} or Chapter II, Section 5 of \cite{So08}.
\end{remark}


\section{Preliminaries}


\subsection{Commutation and null forms}

In preparation for the energy estimates, we need to consider
the commutation properties of the vector fields $\Gamma$ with respect
to the nonlinear terms.  It is necessary to verify that the
null structure is preserved upon differentiation, in some sense.

\begin{lemma}
\label{nd}
Let $u$ be a solution of \eqref{pde}.
Assume that the null condition \eqref{nc} holds for the nonlinearity
in \eqref{non}.
Then for any $\al\in\N^8$,
$$
\Box_\gm \Gamma^\al u = \sum_{\be+\ga+\mu=\al}N_\mu^\al(\Gamma^\be u,\Gamma^\ga u)+\sum_{|\be|\le|\al|-1}\(r_0 \nabla^2 \Gamma^\be u+r_1 \nabla \Gamma^\be u\right),
$$
where each $N_\mu^\al$ is a quadratic nonlinearity of the form
\eqref{non} satisfying \eqref{nc}, and $r_m$ with $m\in\N$, which may change from line to line, denote functions such that
$$|\nabla^\al r_m|\le C_\al \delta \<r\>^{-\rho-m-|\al|}\quad\mbox{for any}\quad\al\in\N^3\ .$$
Moreover, if $|\mu|=0$, then $N_\mu^\al=N$.
\end{lemma}

\begin{prf}
We will give the proof by induction. It is clear that the result is true for $|\al|=0$. Now assume that it is true for any $\al$ with $|\al|=m$. Given $\al_0$ with $|\al_0|=m+1$, we can find some $j$ and $\al$ with $|\al|=m$ and $\Ga^{\al_0}=\Ga_j\Ga^\al$.

To proceed,
we define
$$
[\Gamma,N](u,v)=\Gamma N(u,v) - N(\Gamma u,v)-N(u,\Gamma v),
$$
which is a quadratic nonlinearity of the form \eqref{non}. Moreover, by the proof of Lemma 4.1 in \cite{ST}, we know that $[\Gamma,N]$ is of null form for each $\Gamma$.

By \eqref{H1}, we have
$\Delta_\gm =\Delta+r_0\nabla^2+r_1\nabla$.
Recall also that
$$[\pt^2-c^2_I \Delta, \Ga_j]=2\delta_{j7}(\pt^2-c^2_I \Delta)\ ,$$
then, the term  $(\Box_\gm \Ga^{\al_0} u)^I$  could be calculated as follows
\begin{eqnarray*}
&&(\pt^2-c_I^2 \Delta_\gm) \Ga_j\Ga^\al u^I\\
  &=&[\pt^2-c_I^2 \Delta_\gm,\Ga_j] \Ga^\al u^I+\Ga_j  (\pt^2-c_I^2 \Delta_\gm)\Ga^\al u^I\\
  &=&2\delta_{j7}(\pt^2-c_I^2 \Delta) \Ga^\al u^I-c_I^2[r_0\nabla^2+r_1\nabla,
  \Ga_j]\Ga^\al u^I\\
  &&+ \sum_{\be+\ga+\mu=\al}\Ga_j N_\mu^{\al,I}(\Gamma^\be u,\Gamma^\ga u)+\sum_{|\be|\le|\al|-1}\Ga_j\(r_0 \nabla^2 \Gamma^\be u^I+r_1 \nabla \Gamma^\be u^I\right)\\
  &=&2\delta_{j7}(\Box_\gm \Ga^\al u)^I+\sum_{|\be|\le|\al|}\(r_0 \nabla^2 \Gamma^\be u^I+r_1 \nabla \Gamma^\be u^I\right)\\
  &&+ \sum_{\be+\ga+\mu=\al}\left\{[\Ga_j, N_\mu^{\al,I}](\Gamma^\be u,\Gamma^\ga u)+N_\mu^{\al,I}(\Ga_j\Gamma^\be u,\Gamma^\ga u)+N_\mu^{\al,I}(\Gamma^\be u, \Ga_j\Gamma^\ga u)\right\}
\end{eqnarray*}
which is of the required form,   by the induction assumption. This completes the proof.
\end{prf}

\subsection{Estimates for null forms}
The utility of the null condition is captured in the next lemma, where we get some additional decay in nonlinearities with the null structure \eqref{nc}.

\begin{lemma}
\label{ptwnc}
Suppose that the nonlinear form $N(u,v)$
defined in \eqref{non} satisfies the null condition \eqref{nc}.
For any
$u$, $v$, $w\in C^2([0,T]\times\R^3)$ and $r\ge C_I t/2$, we have
at any point $(t,x)\in [0,T]\times\R^3$
\begin{subequations}
\begin{multline}
\label{ptwnc1}
|Q^{I,\al\be\ga}_{II}\partial_\al u
\partial_\be\partial_\ga v|\\
\le\frac{C}{\langle r\rangle}\Big[|\Gamma u||\partial^2v|
+|\partial u||\pa\Ga v|
+|\partial u||\pa v|
+\langle c_I t-r\rangle|\pa u|
|\partial^2v| \Big]\ ,
\end{multline}
\begin{multline}
\label{ptwnc2}
|Q^{I,\al\be\ga}_{II}\pa_\al u
\partial_\be v \partial_\ga w|\\
\le \frac{C}{\langle r\rangle}\Big[
|\Gamma u| |\partial v| |\partial w|
+ |\partial u| |\Gamma v| |\partial w|
+ |\partial u| |\partial v| |\Gamma w|
+ \langle c_I t-r\rangle |\pa u| |\partial v| |\partial w|\Big]\ ,
\end{multline}
and
\begin{multline}
\label{ptwnc3}
|S^{I,\al\be}_{II}\pa_\al u
\partial_\be v |
\le \frac{C}{\langle r\rangle}\Big[
|\Gamma u| |\partial v|
+ |\partial u| |\Gamma v|
+ \langle c_I t-r\rangle |\partial u| |\pa v| \Big].
\end{multline}
\end{subequations}
\end{lemma}
\begin{prf}
The inequalities \eqref{ptwnc1}-\eqref{ptwnc2} are exactly Lemma 5.1 of \cite{ST}. The proof of
\eqref{ptwnc3} is similar. See also Lemma 5.4 of Chapter II in \cite{So08}.
\end{prf}

\subsection{Sobolev-type inequalities}
\label{sob}

The following Sobolev inequalities do not involve
the Lorentz boost operators.  The weight $\langle c_I t-r\rangle$
compensates for this.  We use the notation defined in
\eqref{ennorm}, \eqref{wnorm}.

\begin{lemma}
We have the following inequalities for smooth functions $u:\R^{3+1}_+\rightarrow \R^M$,
\begin{alignat}{2}
\label{ellinf4}
&\langle r \rangle^{1/2}| u(t,x)|\le C E_2^{1/2}(u(t)),
\\
\label{ellinf5}
&\langle r \rangle |\partial  u(t,x)|\le C E_3^{1/2}(u(t)),
\\
\label{ellinf6}
&\langle r\rangle\langle c_I t -r \rangle^{1/2}
|\partial u^I(t,x)|
\le C\Big[E_3^{1/2}(u(t))+\mathcal{X}_3(u(t))\Big],\\
\label{ellinf7}
&\langle r\rangle\langle c_I t -r \rangle
| \partial^2 u^I(t,x)|
\le C\mathcal{X}_4(u(t)),
\\
\label{ellinf8}
&\langle r\rangle^{1/2}\langle c_I t -r \rangle
|\partial u^I(t,x)|
\le C\Big[E_2^{1/2}(u(t))+\mathcal{X}_3(u(t))\Big].
\end{alignat}
\end{lemma}

See \cite{KlSid} and Proposition 3.3 in  \cite{Sid2} for the proof of \eqref{ellinf4}-\eqref{ellinf7}. The inequality \eqref{ellinf8} is just (4.22) of \cite{Hi04}. See also Lemma 5.2 of Chapter II in \cite{So08}.

\subsection{Local energy estimates}
One of the main extra steps in our proof is to exploit the local energy estimates (also known as Morawetz estimates, KSS estimates), to handle the extra terms arising from the non-flat metric.

\begin{lemma}\label{KSS}
Let $f_0=(r/(1+r))^{2\mu}$, $f_k=r/(r+2^{k})$ with $k\ge 1$ and $\mu\in (0,1/2)$, and $u$ be the solution to the equation $(\pt^2-c_I^2 \Delta +h^{I,\al\be}(t,x)\pa_\al\pa_\be)u^I=F^I$ in $[0,T]\times \R^n$ with $h^{I,\al\be}=h^{I,\be\al}$, 
$|h^{I,00}|\le 1/2$, $\sum_{1\le i,j\le n} |h^{I,ij}|\le c_I^2 /2$ for any $I$
 and $n\ge 3$, then there exists a constant $C>0$, depending only on the dimension $n$, such that
\begin{eqnarray*}
&&\sup_{0\le t\le T} E_1(u(t))+LE_1(T)+(\log(2+T))^{-1}KSS_1(T)\\
\\
&\le &C E_1(u(0))+
C\int_0^T \int_{\R^n}
\left\{
   \left( |\pa h|+\frac{|h|}{r^{1-2\mu}\<r\>^{2\mu}}   \right)
   |\pa u|
    \left( |\pa u|+\frac{|u|}{r}\right)
\right\}
dx dt\\
&&+C\left|\sum_{I=1}^M\int_0^T \int_{\R^n} F^I \pt u^I  dx dt\right|
+C\sup_{k\ge 0}\left|\sum_{I=1}^M\int_0^T \int_{\R^n}
	f_k
		\left(\pa_r u^I+\frac{n-1}{2r}u^I\right)
	F^I dx dt\right|
\ .
\end{eqnarray*}
\end{lemma}

See Section 2 of \cite{HWY10} for the proof (see also \cite{MeSo06, HWY11} and references therein).

\subsection{Weighted decay estimates}
One of the main steps is to control the weighted norm $\mathcal{X}_\od(u(t))$.
This is accomplished in this subsection by a type of bootstrap argument, similar to that in \cite{ST}.
\begin{lemma}[Klainerman-Sideris estimate]
\label{decay1}
Let $\de\ll 1$ and $\rho\ge 1$.  Then
\begin{equation}
\label{locen}
\mathcal{X}_\od(u(t))\le C\left[E_{\od}^{1/2}(u(t))+
\sum_{|\al|\le\od-2}\|(t+r)\,\Box_\gm\Gamma^\al u(t)\|_{L^2}\right]
\end{equation}
for any $u\in \dot{H}_\Gamma^\od(T)$.
\end{lemma}

\begin{prf}
The same estimate is known to be true for the standard D'Alembertian $\Box$ instead of $\Box_\gm$, see Lemma 7.1 in \cite{ST},  \cite{KlSid} and Lemma 5.3 of Chapter II in \cite{So08}. To complete the proof, we need only to control the norm involving $\Box$ by that of $\Box_\gm$ and good terms, as follows. For any $\al\in\N^8$ with $|\al|\le \od-2$,
\begin{eqnarray*}
&&	\|(t+r)\,\Box\Gamma^\al u(t)\|_{L^2}-\|(t+r)\,\Box_\gm\Gamma^\al u(t)\|_{L^2}\\
&\le&
	\|(t+r)\, r_0\nabla^2 \Gamma^\al u(t)\|_{L^2}+\|(t+r)\, r_1\nabla \Gamma^\al u(t)\|_{L^2}\\
&\le&
	\sum_{I}\left(\|\<c_I t-r\>\, r_{-1}\nabla^2 \Gamma^\al u^I(t)\|_{L^2}+\|\<c_I t-r\>\, r_0\nabla \Gamma^\al u^I(t)\|_{L^2}\right)\\
&\le&
	C \de \X_\od(u(t))+C \de \|\nabla (\<c_I t-r\>\nabla \Gamma^\al u^I(t))\|_{L^2}\\
&\le&
C \de \X_\od(u(t))+C \de \|\nabla \Gamma^\al u(t)\|_{L^2}\\
&\le&
C \de \X_\od(u(t))+C \de E_\od^{1/2}(u(t))\ ,
\end{eqnarray*}
where we have used the elementary inequality \beeq\label{elem}t+r\le C\<c_I t-r\>\<r\>\eneq in the second inequality, Hardy's inequality and $\rho\ge 1$ in the third inequality, and the fact that $|\nabla \<c_I t-r\>|\le C$ in the fourth inequality.
\end{prf}

Now we assume that $u$ solves the nonlinear equation \eqref{pde}.

\begin{lemma}
\label{lem2}
Let $\rho\ge 1$ and $u\in \dot{H}_\Gamma^\od(T)$ be a solution of \eqref{pde}.
Define ${\od}'=\left[\frac{\od-1}{2}\right]+3.$
Then  for all $|\al|\le\od-2$,
\begin{multline}
\label{w2}
\|(t+r)\Box_\gm \Gamma^\al u(t)\|_{L^2}
\le C[\X_{\od'}(u(t))
E_{\od-1}^{1/2}(u(t))
+\X_{\od}(u(t))
E_{{\od}'}^{1/2}(u(t))]\\
+CE_{\od'}^{1/2}(u(t))E_{\od-1}^{1/2}(u(t))
+C\de \X_{\od-1}(u(t))+C\de E_{\od-2}^{1/2}(u(t)).
\end{multline}
\end{lemma}
\begin{prf}
There are similar estimates for the case of $\Box$ in Lemma 7.2 of \cite{ST} and Lemma 5.2 of \cite{Hi04}.
In view of Lemma \ref{nd}, we need to control the terms of the form
$$\|(t+r) \pa \Gamma^\be u^I (t) \pa^2 \Ga^\ga u^J(t) \|_{L^2}\ ,$$
$$\|(t+r) \pa \Gamma^\be u^I(t) \pa \Ga^\ga u^J(t) \|_{L^2}\ ,$$
with $|\be|+|\ga|\le|\al|\le\od-2$ and
$$\|(t+r) r_0 \nabla^2 \Ga^\be u(t) \|^2_{L^2} +\|(t+r) r_1 \nabla \Gamma^\be u(t)  \|^2_{L^2}$$
with $|\be|\le |\al|-1\le\od-3$.
For the first set of terms, we separate two cases: either $|\be|\le \od'-3$ or $|\ga|\le \od'-4$. In the case of $|\be|\le \od'-3$, using \eqref{elem} and \eqref{ellinf5}, we have
\begin{eqnarray*}
  \|(t+r) \pa \Gamma^\be u^I(t) \pa^2 \Ga^\ga u^J(t) \|_{L^2}
  &\le & C \|\<r\> \pa \Gamma^\be u^I(t)\|_{L^\infty} \|\<c_J t-r\>\pa^2 \Ga^\ga u^J(t) \|_{L^2}\\
  &\le &C E_{\od'}^{1/2}(u(t)) \X_\od (u(t)) \ .
\end{eqnarray*}
In the second case, using \eqref{elem} and \eqref{ellinf7}, we get
\begin{eqnarray*}
  \|(t+r) \pa \Gamma^\be u^I (t) \pa^2 \Ga^\ga u^J (t) \|_{L^2}
  &\le & C \|\pa \Gamma^\be u^I (t)\|_{L^2} \|\<r\> \<c_J t-r\>\pa^2 \Ga^\ga u^J(t) \|_{L^\infty}\\
  &\le &C E_{\od-1}^{1/2}(u(t)) \X_{\od'} (u(t)) \ .
\end{eqnarray*}

Turning to the second set of terms, without loss of generality, we can assume $|\be|\le |\ga|$ (and so $|\be|\le \od'-3$). Then,
using \eqref{ellinf8} and \eqref{ellinf5}, we get
\begin{eqnarray*}
&&  \|(t+r) \pa \Gamma^\be u^I (t) \pa \Ga^\ga u^J(t) \|_{L^2}\\
  &\le & C (\|\<c_I t-r\>\pa \Gamma^\be u^I(t)\|_{L^\infty}+\|\<r\>\pa \Gamma^\be u^I(t)\|_{L^\infty}) \|\pa \Ga^\ga u^J(t) \|_{L^2}\\
  &\le &C ( E_{\od'-1}^{1/2}(u(t))+ \X_{\od'} (u(t)) ) E_{\od-1}^{1/2}(u(t))+C E_{\od'}^{1/2}(u(t)) E_{\od-1}^{1/2} (u(t)) \ .
\end{eqnarray*}
Since $\rho\ge 1$ and $|\be|\le |\al|-1\le \od-3$, we see that
\begin{eqnarray*}
&&  \|(t+r) r_0 \nabla^2 \Ga^\be u^I(t) \|_{L^2} +\|(t+r) r_1 \nabla \Gamma^\be u^I(t)  \|_{L^2}\\
&\le& C\|\<c_I t-r\> r_{-1} \nabla^2 \Ga^\be u^I(t) \|_{L^2} +\|\<c_I t-r\> r_0 \nabla \Gamma^\be u^I(t)  \|_{L^2}\\
&\le& C\de \X_{\od-1}( u(t))+C\de \|\nabla (\<c_I t-r\>  \nabla \Gamma^\be u(t)  ) \|_{L^2}\\
&\le& C\de \X_{\od-1}( u(t))+C\de E_{\od-2}^{1/2}(u(t))\ ,\end{eqnarray*}
where we have used the Hardy inequality in the second inequality. This completes the proof.
\end{prf}

The next result gains control of the weighted norm $\X$ by the energy.
We distinguish two different energies, the lower order of which will
remain small.  In Section \ref{sec-4}, we will
allow the energy of higher order to grow in time.
\begin{lemma}
\label{we}
Let $u\in \dot{H}_\Gamma^\od(T)$, $\od\ge8$, be a solution of \eqref{pde}.
Define $\eta=\od-2$, and assume that
\[
\de\ll 1,\ \eps_0\equiv\sup_{0\le t < T} E_{\eta}^{1/2}(u(t))\ll 1,\ \rho\ge 1.
\]
Then for $0\le t< T$,
\begin{subequations}
\begin{equation}
\label{we1}
\X_\eta(u(t))\le C E_{\eta}^{1/2}(u(t))
\end{equation}
and
\begin{equation}
\label{we2}
\X_\od(u(t))\le C E_{\od}^{1/2}(u(t)).
\end{equation}
\end{subequations}
\end{lemma}
\begin{prf}
Let $\eta' = \left[\frac{\eta-1}{2}\right]+3$, $\eta=\od-2$.
Since $\eta\ge 6$, we have $\eta'\le\eta$.  Thus, by Lemmas
\ref{decay1} and \ref{lem2}, we get from our assumption
\begin{eqnarray*}
 \X_\eta(u(t))
 &\le& C [E_{\eta}^{1/2}(u(t)) + E_{\eta'}^{1/2}\X_\eta+\X_{\eta'}E^{1/2}_{\eta-1}+E_{\eta'}^{1/2}E^{1/2}_{\eta-1}]\\
 &&+C\de [\X_{\eta-1}(u(t))+ E_{\eta-2}^{1/2}(u(t))]\\
 &\le& C [E_{\eta}^{1/2}(u(t)) + (\ep_0+\de)\X_\eta(u(t))]
 \end{eqnarray*}
Thus, if $\eps_0$ and $\de$ are small enough, the bound \eqref{we1} results.

Again since $\od\ge8$, we have $\od'=\left[\frac{\od-1}{2}\right]+3
\le \eta =\od-2$.   From Lemmas \ref{decay1} and \ref{lem2}, we have
\begin{eqnarray*}
 \X_\od(u(t))
 &\le& C [E_{\od}^{1/2}(u(t)) +\X_{\eta} E_{\od-1}^{1/2}
+\X_{\od} E_{\eta}^{1/2}+E_{\eta}^{1/2}E_{\od-1}^{1/2}]\\
& & + C\de [\X_{\od-1}(u(t)) + E_{\od-2}^{1/2}(u(t))]\\
 &\le& C (1+\X_{\eta}+E_{\eta}^{1/2}+\de) E_{\od}^{1/2}(u(t))
+C( E_{\eta}^{1/2}+\de)\X_{\od}(u(t))
\end{eqnarray*}
If we apply \eqref{we1}, $\vep_0\ll 1$ and $\de\ll 1$, then
\[
\X_\od(u(t))\le C E_{\od}^{1/2}(u(t))
+C(\eps_0+\de)\mathcal{X}_{\od}(u(t)),
\]
from which \eqref{we2} follows.
\end{prf}

\section{Energy and local energy estimates}\label{sec-energy}
Assume that $u(t)\in \dot{H}^\od_\Gamma(T)$ is a local
solution of the initial value problem for \eqref{pde} (in which we need the symmetry condition \eqref{symmetric}).
Our task will be to show that $E_\od(u(t))$ remains finite
for all $t\ge0$.  To do so, we will derive a pair
of coupled integral inequalities for
$E_\od(u(t))+LE_\od(t)$ and $E_{\eta}(u(t))+LE_\eta(t)$, with $\eta=\od-2$.
If \eqref{smallness} holds with $2C_1 \eps\le \ep_0$ occurred in Lemma \ref{we},
then $E_\eta^{1/2}(u(0))<\eps$ and $E_\eta^{1/2}(u(t))<2 C_1\eps$ for certain small interval $t\in [0,T]$.
Define
$$T_0=\sup\{T: E_\eta^{1/2}(u(t))\le 2 C_1\eps, t\in [0,T]\}\ .$$
Here the constants $C_0$, $C_1\ge 1$ will be determined later (see Subsection \ref{sec-conclusion}).
All of the following computations will be valid on the interval $[0,T_0)$.

To complete the proof of Theorem \ref{thm}, we need to obtain estimates for $E_m(u(t))$, as well as $\I_m(u(t))$ with $m=\eta, \od$. Instead of giving the estimates for $\I_m(u(t))$ directly, we will give the estimates of $LE_m(u(t))$.

Since we have $$ (\Box_\gm)^I = \pt^2-c_I^2 \Delta - c_I^2 (g^{ij}-\delta^{ij})\pa_i\pa_j+r_1 \nabla \ ,$$
it is easy to see that Lemma \ref{KSS} applies to $(\Box)^I-h^{I,ij}\pa_i\pa_j$ with $h^{I,ij}=c_I^2 (g^{ij}-\delta^{ij})$.
By Lemma \ref{nd},  we see that
\begin{eqnarray*}
&&(\pt^2-c_I^2 \Delta - c_I^2 (g^{ij}-\delta^{ij})\pa_i\pa_j)\Ga^\al u^I\\
&=&(\Box_\gm \Ga^\al u)^I-r_1\nabla \Ga^\al  u^I\\
&=&\sum_{\be+\ga+\mu=\al}N_\mu^{\al,I}(\Gamma^\be u,\Gamma^\ga u)+\left[\sum_{|\be|\le|\al|-1}\(r_0 \nabla^2 \Ga^\be u^I+r_1 \nabla \Ga^\be u^I\right)-r_1\nabla \Ga^\al u^I\right]\\
&=&F_\al^I+G_\al^I\ .
\end{eqnarray*}
Since $\rho>1$ and $\mu\in (0,1/2)$, if we assume \beeq\label{120630-1}2\mu'\le \rho-1+2\mu\ ,\eneq
then Lemma \ref{KSS} tells us that
\begin{eqnarray}
  &&\sup_{t\in [0,T]}E_m(u(t))+LE_m(T)\label{KSS2}\\
  &\le & C E_m(u(0))+ C\sum_{|\al|\le m-1} \int_0^T \int_{\R^3}
\left(|\pa \Ga^\al u|+\frac{1}{r}  |\Ga^\al u | \right)
	|G_\al | dx dt\nonumber\\
  &&+C\de
\sum_{|\al|\le m-1} \int_0^T \int_{\R^3}
\left\{
   r^{-1+2\mu}\<r\>^{-2\mu-\rho}    |\pa \Ga^\al u|
    \left( |\pa \Ga^\al u|+\frac{|\Ga^\al u|}{r}\right)
\right\} dx dt \nonumber\\
&&+C\sum_{|\al|\le m-1} \left|\sum_{I=1}^M\int_0^T \int_{\R^3} F^I_\al \pt \Ga^\al u^I  dx dt\right| \nonumber\\
&&
+C\sum_{|\al|\le m-1} \sup_{k\ge 0}\left|\sum_{I=1}^M\int_0^T \int_{\R^3}
	f_k
		\left(\pa_r \Ga^\al u^I+\frac{1}{r}  \Ga^\al u^I\right)
	F^I_\al dx dt\right|\nonumber\\
&\le &C E_m(u(0))+C \de LE_m(T)+
C\sum_{|\al|\le m-1}\left|\int_0^T C_1^{\al}(t)dt\right|+\sup_{k\ge 0}\left|\int_0^T C_{2,k}^{\al}(t)dt\right|\ .\nonumber  \end{eqnarray}
Here, we have introduced the notations
\beeq\label{0609-1}C_1^{\al}=\int_{\R^3} \sum_{I=1}^M F^I_\al \pt \Ga^\al u^I dx
\ ,
\eneq
\beeq\label{0609-2}C_{2,k}^{\al}(t)=
\int_{\R^3}\sum_{I=1}^M
	f_k
		\left(\pa_r \Ga^\al u^I+\frac{1}{r}  \Ga^\al u^I\right)
	F^I_\al dx \ ,
\eneq
where
$F^I_\al=\sum_{\be+\ga+\mu=\al}N_\mu^{\al,I}(\Gamma^\be u,\Gamma^\ga u)$
and $N_\mu^{\al,I}$ denotes the $I$-th component of $N_\mu^{\al}=Q_\mu^{\al}+S_\mu^{\al}$ (see \eqref{non}).

\subsection{Estimate for $C_1^{\al}$}
From \eqref{0609-1} we know
$$C_1^{\al}=\int_{\R^3} \sum_{I=1}^M F^I_\al \pt \Ga^\al u^I dx=
\sum_{I=1}^M\sum_{\be+\ga+\mu=\al} \int_{\R^3}  N_\mu^{\al,I}(\Gamma^\be u,\Gamma^\ga u)  \pt \Ga^\al u^I dx\ .
$$ Among all of the terms, the cases $|\al|=|\ga|=m-1$ for $Q_\mu^\al$ are quasilinear terms and we want to use the symmetry condition \eqref{symmetric} to absorb such terms. When $|\al|=|\ga|=m-1$, we have $Q^\al_\mu=Q$, $\ga=\al$. Then an integration by parts argument (see \cite{ST}) will yield
\begin{eqnarray*}
&& \sum_{I=1}^M\int_{\R^3}  Q^I(u,\Gamma^\al u) \pt \Ga^\al u^I(t) dx\\
&=&\sum_{I=1}^M Q^{I,\be\mu \ga}_{JK}  \int_{\R^3}\left[\pa_\ga\left( \pa_\be u^J \pa_{\mu} \Ga^\al u^K \pt \Ga^\al u^I\right) \right.\\
&&\left.-  \pa_\ga \pa_\be u^J \pa_{\mu} \Ga^\al u^K \pt \Ga^\al u^I -  \pa_\be u^J \pa_{\mu} \Ga^\al u^K  \pa_\ga \pt \Ga^\al u^I \right]dx
\\
&=& \sum_{I=1}^M Q^{I,\be\mu 0}_{JK} \pt \int_{\R^3} \pa_\be u^J \pa_{\mu} \Ga^\al u^K \pt \Ga^\al u^I dx\\
&&
-\sum_{I=1}^M Q^{I,\be\mu\ga}_{JK}  \int_{\R^3} \left[\pa_\ga \pa_\be u^J \pa_{\mu} \Ga^\al u^K \pt \Ga^\al u^I
+\frac{1}{2}\pa_\be u^J \pt\left(\pa_{\mu} \Ga^\al u^K  \pa_\ga  \Ga^\al u^I \right)\right] dx\\
&=& \sum_{I=1}^M(Q^{I,\be\mu 0}_{JK}\delta_0^\ga-\frac 12 Q^{I,\be\mu\ga}_{JK}) \pt \int_{\R^3} \pa_\be u^J \pa_{\mu} \Ga^\al u^K
\pa_\ga\Ga^\al u^I  dx\\
& &- \sum_{I=1}^MQ^{I,\be\mu\ga}_{JK}  \int_{\R^3} \left[\pa_\ga \pa_\be u^J \pa_{\mu} \Ga^\al u^K \pt \Ga^\al u^I -\frac{1}{2}  \pa_\be\pt u^J \pa_{\mu} \Ga^\al u^K  \pa_\ga  \Ga^\al u^I \right]dx
\end{eqnarray*}
where we have used the symmetry of the equation \eqref{symmetric} in the second step. By introducing the notation \beeq\label{newEnergy1}\tilde E_{\al}(t)=
\sum_{I=1}^M (Q^{I,\be\mu 0}_{JK}\delta_0^\ga-\frac 12 Q^{I,\be\mu\ga}_{JK}) \int_{\R^3} \pa_\be u^J \pa_{\mu} \Ga^\al u^K
\pa_\ga\Ga^\al u^I  dx\ ,\eneq
we see that
\begin{eqnarray}
&& \sum_{|\al|\le m-1} \left|\int_0^T C_1^{\al}(t)dt\right| \label{localenergy1}\\
&\le &
\sum_{I=1}^M \sum_{|\al|\le m-1, \be+\ga+\mu=\al} \left|\int_0^T\int_{\R^3} S_\mu^{\al,I}(\Gamma^\be u,\Gamma^\ga u) \pt \Ga^\al u^I(t) dx dt\right|
\nonumber\\
&&+ \sum_{I=1}^M \sum_{|\al|\le m-1, \be+\ga+\mu=\al, |\ga|<m-1} \left|\int_0^T\int_{\R^3} Q_\mu^{\al,I}(\Gamma^\be u,\Gamma^\ga u) \pt \Ga^\al u^I(t) dx dt\right|
\nonumber\\
&&+ \sum_{I=1}^M \sum_{|\al|=m-1}\left|\int_0^T \int_{\R^3}Q^{I,\be\mu\ga}_{JK}   \pa_\ga \pa_\be u^J \pa_{\mu} \Ga^\al u^K \pt \Ga^\al u^I dxdt \right|\nonumber\\
&&+  \sum_{I=1}^M  \sum_{|\al|=m-1}\left|\int_0^T \int_{\R^3} Q^{I,\be\mu\ga}_{JK} \pa_\be\pt u^J \pa_{\mu} \Ga^\al u^K  \pa_\ga  \Ga^\al u^I dx dt \right| \nonumber\\
&&+ \sum_{|\al|=m-1} |\tilde E_{\al}(T)-\tilde E_{\al}(0)|\ .
\nonumber
\end{eqnarray}

\subsection{Estimate for $C_{2,k}^{\al}$}
Let $L_k=f_k 	(\pa_r +\frac{1}{r} )=h_k^i(x)\pa_i+ h_k(x)$ with $h_k(x)=f_k/r$ and $h_k^i(x)=f_k x^i/r$,
we know from \eqref{0609-2} that
$$C_{2,k}^{\al}=\int_{\R^3}
	\sum_{I=1}^M 	F_\al^I L_k \Ga^\al u^I  dx=\sum_{I=1}^M
\sum_{\be+\ga+\mu=\al} \int_{\R^3}  N_\mu^{I,\al}(\Gamma^\be u,\Gamma^\ga u)  L_k \Ga^\al u^I dx\ .
$$
As for $C_1^\al$, for the case $|\al|=|\ga|=m-1$, we have $Q^{\al}_\mu=Q$,  $\ga=\al$, and
\begin{eqnarray}
&& \sum_{I=1}^M \int_{\R^3} Q^I(u,\Gamma^\al u) L_k \Ga^\al u^I(t) dx\\
  &=& \sum_{I=1}^M \int_{\R^3} Q^{I,\be\mu\ga}_{JK}\pa_\be u^J \pa_{\mu\ga} \Ga^\al u^K  L_k \Ga^\al u^I dx \nonumber\\
&=&\sum_{I=1}^M  Q^{I,\be\mu 0}_{JK} \pt \int_{\R^3} \pa_\be u^J \pa_{\mu} \Ga^\al u^K  L_k \Ga^\al u^I dx\nonumber\\
&&- \sum_{I=1}^M Q^{I,\be\mu\ga}_{JK}  \int_{\R^3} \pa_\ga \pa_\be u^J \pa_{\mu} \Ga^\al u^K  L_k \Ga^\al u^I dx\nonumber\\
&&- \sum_{I=1}^M Q^{I,\be\mu\ga}_{JK}  \int_{\R^3} \pa_\be u^J \pa_{\mu} \Ga^\al u^K  \pa_\ga  L_k \Ga^\al u^I dx\ .\label{605-1}
\end{eqnarray}
The last term \eqref{605-1} in the above identity can be rewritten as follows,
\begin{eqnarray*}
&&- \sum_{I=1}^M Q^{I,\be\mu\ga}_{JK}\int_{\R^3} \pa_\be u^J \pa_{\mu} \Ga^\al u^K  \left[ (\pa_\ga h^i_k)\pa_i  +h^i_k \pa_\ga \pa_i +  (\pa_\ga h_k)+  h_k \pa_\ga\right]\Ga^\al u^I dx\\
&=&-\sum_{I=1}^M  Q^{I,\be\mu\ga}_{JK}\int_{\R^3} \pa_\be u^J \pa_{\mu} \Ga^\al u^K  \left[ (\pa_\ga h^i_k)\pa_i  +(\pa_\ga h_k)+ h_k \pa_\ga\right]\Ga^\al u^I dx\nonumber\\
& &- \frac 12 \sum_{I=1}^M Q^{I,\be\mu\ga}_{JK} \int_{\R^3}h^i_k  \pa_\be u^J \pa_i ( \pa_{\mu} \Ga^\al u^K   \pa_\ga  \Ga^\al u^I) dx\nonumber\\
&=&\frac 12\sum_{I=1}^M   Q^{I,\be\mu\ga}_{JK} \int_{\R^3}
h^i_k  \pa_\be  \pa_i u^J \pa_{\mu} \Ga^\al u^K   \pa_\ga \Ga^\al u^I dx\\
&& - \sum_{I=1}^M Q^{I,\be\mu\ga}_{JK}\int_{\R^3} \pa_\be u^J \pa_{\mu} \Ga^\al u^K  \left[ (\pa_\ga h^i_k)\pa_i  + \left(h_k-\frac{1}2 \pa_i h^i_k\right)  \pa_\ga+\pa_\ga h_k\right]\Ga^\al u^I dx
\end{eqnarray*}
where we have used the symmetric condition \eqref{symmetric}.
 It is easy to check that, by setting, say, $\mu=1/4$, there is a uniform constant $C>0$ which is independent of $k\ge 0$ such that
$$|\nabla h_k^i|+|h_k|\le C r^{-1/2}\<r\>^{-1/2},\ |\nabla h_k|\le C r^{-3/2}\<r\>^{-1/2}.$$
Moreover, by \eqref{ellinf5}, we have
$$|\pa u|\le C \<r\>^{-1} E_3^{1/2}(u(t))\le CC_1\ep \<r\>^{-1}\ .$$
Then,  the last term in the above expression is controlled by $CC_1\ep \I_m(t)$ with \beeq\label{120630-2}\mu=1/4<\mu'\le 3/4\ .\eneq

In summary, by introducing the notation \beeq\label{newEnergy2}
 \tilde E_{\al,k}=\sum_{I=1}^M
Q^{I,\be\mu 0}_{JK} \int_{\R^3} \pa_\be u^J \pa_{\mu} \Ga^\al u^K  L_k \Ga^\al u^I dx
, \eneq
 we have proven the following
\begin{eqnarray}
&&\quad\quad \sum_{|\al|\le m-1} \left|\int_0^T C_{2,k}^{\al}(t)dt\right| \label{localenergy2}\\
&\le &\sum_{I=1}^M
 \sum_{|\al|\le m-1, \be+\ga+\mu=\al} \left|\int_0^T\int_{\R^3} S_\mu^{\al,I}(\Gamma^\be u,\Gamma^\ga u) L_k \Ga^\al u^I(t) dx dt\right|
\nonumber\\
&&+ \sum_{|\al|\le m-1, \be+\ga+\mu=\al, |\ga|<m-1} \left|\int_0^T\int_{\R^3} Q_\mu^{\al,I}(\Gamma^\be u,\Gamma^\ga u) L_k \Ga^\al u^I(t) dx dt\right|
\nonumber\\
&&+ \sum_{I=1}^M \sum_{|\al|=m-1}\left|\int_0^T \int_{\R^3}Q^{I,\be\mu\ga}_{JK}   \pa_\ga \pa_\be u^J \pa_{\mu} \Ga^\al u^K L_k \Ga^\al u^I dxdt \right|\nonumber\\
&&+ \sum_{I=1}^M   \sum_{|\al|=m-1}\left|\int_0^T \int_{\R^3} Q^{I,\be\mu\ga}_{JK} h^i_k\pa_\be\pa_i u^J \pa_{\mu} \Ga^\al u^K  \pa_\ga  \Ga^\al u^I dx dt \right| \nonumber\\
&&+ \sum_{|\al|=m-1} | \tilde E_{\al,k}(T)- \tilde E_{\al,k}(0)|+CC_1\ep LE_m(T)\ .
\nonumber
\end{eqnarray}

\section{Proof of the energy inequality}\label{sec-4}
In this section, we will complete the proof of Theorem \ref{thm}.

At first, we claim that, with $\od\ge 9$ and $\eta=\od-2$, by \eqref{KSS2}, \eqref{localenergy1} and \eqref{localenergy2}, we can deduce that
\begin{eqnarray}
&&LE_\od(T)+\sup_{0\le t\le T}E_\od(u(t))\label{en2}\\
&\le& C_3 E_\od(u(0))+C_3 \int_0^T \langle t \rangle^{-1} {E}_\eta^{1/2}(u(t))
{E}_\od(u(t)) dt\nonumber\\
&&+C_3C_1(\de+\ep) LE_\od(T)+C_3\sum_{|\al|=\od-1} \sup_{k\ge 0}|E_{2,k}^{\al}(T)-E_{2,k}^{\al}(0)|\nonumber\\
&&+C_3\sum_{|\al|=\od-1} |E_{1}^{\al}(T)-E_{1}^{\al}(0)|\ ,
\nonumber
\end{eqnarray}
and
\begin{eqnarray}
&&LE_\eta(T)+\sup_{0\le t\le T}E_\eta(u(t))\label{en3}\\
&\le& C_3 E_\eta(u(0))+C_3 \int_0^T \langle t \rangle^{-1}t^{-1/2} {E}_\od^{1/2}(u(t))
{E}_\eta(u(t)) dt\nonumber\\
&&+C_3C_1(\de+\ep) LE_\eta(T)+C_3\sum_{|\al|=\eta-1} \sup_{k\ge 0}|E_{2,k}^{\al}(T)-E_{2,k}^{\al}(0)|\nonumber\\
&&+C_3\sum_{|\al|=\eta-1} |E_{1}^{\al}(T)-E_{1}^{\al}(0)|\ .
\nonumber
\end{eqnarray}
for some universal constant $C_3\ge 1$.

Comparing the right hand sides in the estimates  \eqref{localenergy1} and \eqref{localenergy2}, we find similar terms involving the integral of
\beeq\label{term1}
D_1^1=Q_\mu^{\al,I}(\Gamma^\be u,\Gamma^\ga u) \pa \Ga^\al u^I(t) ,\  D_1^2=h_k(x) Q_\mu^{\al,I}(\Gamma^\be u,\Gamma^\ga u)  \Ga^\al u^I(t)
\eneq
for $|\al|\le m-1$, $\be+\ga+\mu=\al$ and $|\ga|<m-1$,
\beeq\label{term2}
D_2^1=
Q^{I,\be\mu\ga}_{JK}   \pa_\ga \pa_\be u^J \pa_{\mu} \Ga^\al u^K \pa \Ga^\al u^I,
\
D_2^2=h_k(x)
Q^{I,\be\mu\ga}_{JK}   \pa_\ga \pa_\be u^J \pa_{\mu} \Ga^\al u^K  \Ga^\al u^I
\eneq
for $|\al|=m-1$,
\beeq\label{term3}
D_3^1=
Q^{I,\be\mu\ga}_{JK} \pa_\be\pa u^J \pa_{\mu} \Ga^\al u^K  \pa_\ga  \Ga^\al u^I
,\
D_3^2=0
\eneq
for $|\al|=m-1$, and the semilinear terms
\beeq\label{term4}
D_4^1=S_\mu^{\al,I}(\Gamma^\be u,\Gamma^\ga u) \pa \Ga^\al u^I(t),\
D_4^2=h_k(x) S_\mu^{\al,I}(\Gamma^\be u,\Gamma^\ga u)  \Ga^\al u^I(t)
\eneq
for $|\al|\le m-1$, $\be+\ga+\mu=\al$.

To prove the claim, we need only to control these terms for $m=\od, \eta$ separately.
The estimates of the first components of $D_j^1$, i.e., $(D_j^1)_0$, have been obtained in \cite{ST} and \cite{Hi04}. Here, we give the proof of the bound for the new terms occurred on the right hand side of \eqref{localenergy2}, that is, to give uniform (with respect to $k$) control of
$$I_j=\left|\int_{\R^3} h_k^i (D_j^1)_i dx\right|,  II_j=\left|\int_{\R^3}  D_j^2 dx\right|$$
with $j=1,2,3,4$.
The proof will be similar to the proof of (8.5) and (8.6) in \cite{ST}.

\subsection{Higher energy}
For the first series of estimates we take $m=\od$ in \eqref{term1}-\eqref{term4}. Recall that $|h_k^i|\le 1$, $|h_k(x)|\le r^{-1/2}(1+r)^{-1/2}$ and using Hardy's inequality,
we obtain immediately
\begin{eqnarray}
I_j+II_j
&\le& C
\sum_{I,J=1}^M\sum_{|\al|\le\od-1,\be+\ga\le \al, |\ga|\le\od-2}
\|\partial\Gamma^\be u^I \partial^2\Gamma^\ga u^J\|_{L^2}
\|\partial\Gamma^\al u\|_{L^2}\label{en1}
\\
&&+
C \sum_{I,J=1}^M\sum_{|\al|\le\od-1,\be+\ga\le \al}
\|\partial\Gamma^\be u^I \partial\Gamma^\ga u^J\|_{L^2}
\|\partial\Gamma^\al u\|_{L^2}
.\nonumber
\end{eqnarray}
For the first term on the right-hand side of \eqref{en1},
 we have either $|\be|\le \od'$
or $|\ga|\le \od'-1$, with $\od'=\left[\od/2\right]$.
Note that since $\od\ge 9$, we have $\od'+3\le\od-2=\eta$.
We will also use that $\langle t \rangle
\le C \langle r \rangle \langle c_J t-r \rangle$.

In the first case, we estimate using \eqref{ellinf5} and \eqref{we2}
\begin{align*}
\|\partial\Gamma^\be u^I \partial^2 \Ga^\ga u^J\|_{L^2}
\le & C\langle t \rangle^{-1}
\|\langle r \rangle \partial \Ga^\be  u^I\|_{L^\infty}
\|\langle c_J t-r \rangle \partial^2 \Ga^\ga u^J\|_{L^2}\\
\le & C\langle t \rangle^{-1} E_{|\be|+3}^{1/2}(u(t)) \mathcal{X}_\od(u(t))\\
\le & C\langle t \rangle^{-1} E_\eta^{1/2}(u(t)) E_\od^{1/2}(u(t)).
\end{align*}
In the second case, we use \eqref{ellinf7} and then \eqref{we1}
\begin{align*}
\|\partial \Ga^\be  u^I \partial^2 \Ga^\ga u^J\|_{L^2}
\le & C\langle t \rangle^{-1}
\|\partial \Ga^\be  u\|_{L^2}
\|\langle r \rangle\langle c_J t-r \rangle \partial^2 \Ga^\ga u^J\|_{L^\infty}\\
\le & C\langle t \rangle^{-1}
E_\od^{1/2}(u(t)) \mathcal{X}_{|\ga|+4}(u(t))\\
\le & C\langle t \rangle^{-1}
E_\od^{1/2}(u(t)) \mathcal{X}_{\eta}(u(t))\\
\le & C\langle t \rangle^{-1}
E_\od^{1/2}(u(t)) E_\eta^{1/2}(u(t)).
\end{align*}

For the second term on the right-hand side of \eqref{en1},
as in the proof of Lemma \ref{lem2},
we can use \eqref{ellinf8} and \eqref{ellinf7} to get (assuming $|\be|\le |\ga|$)
\begin{eqnarray*}
&&  \| \pa \Gamma^\be u^I (t) \pa \Ga^\ga u^J(t) \|_{L^2}\\
  &\le & C\<t\>^{-1} (\|\<c_I t-r\>\pa \Gamma^\be u^I(t)\|_{L^\infty}+\|\<r\>\pa \Gamma^\be u^I(t)\|_{L^\infty}) \|\pa \Ga^\ga u^J(t) \|_{L^2}\\
  &\le &C \<t\>^{-1} ( E_{\eta}^{1/2}(u(t))+ \X_{\eta} (u(t)) ) E_{\od}^{1/2}(u(t))+C E_{\eta}^{1/2}(u(t)) E_{\od}^{1/2} (u(t)) \\
&\le& C\<t\>^{-1}  E_{\eta}^{1/2}(u(t)) E_{\od}^{1/2} (u(t))   \ .
\end{eqnarray*}

Going back to \eqref{en1}, we have established the inequality \eqref{en2}.

\subsection{Lower energy}
The second series of energy type estimates with $m=\eta$ will exploit the
null condition.
Let $c_0=\min\{c_I\}$, then the integrals will be subdivided into separate integrals over the
regions $r\le c_0 t/2$ and $r\ge c_0 t/2$.

{\bf Inside the cones}.
On the region $r\le c_0 t/2$, we have that
\begin{eqnarray*}
\sum_j I_j(r\le c_0 t/2) &\le& C \sum_{I,J,K}\sum_{\be+\ga\le\al, |\al|\le \eta-1,  \ga\le \eta-2}
\|\partial \Ga^\be u^I \partial^2 \Ga^\ga u^J \partial\Gamma^\al u^K \|_{L^1(r\le c_0
t/2)}\\
&&+C
\sum_{I,J,K}\sum_{|\al|\le\eta-1,\be+\ga\le \al}
\|\partial\Gamma^\be u^I \partial\Gamma^\ga u^J
\partial\Gamma^\al u^K\|_{L^1(r\le c_0 t/2)}.
\end{eqnarray*}
Since $r\le c_0 t/2$, we have
 $ \langle t \rangle\le C\langle c_I t-r\rangle
$ for any $I$.  Recall also $|\be|+3\le\od$, $|\ga|+2\le\eta$, and $|\al|+1\le\eta$, thus, using
\eqref{ellinf6} and Lemma \ref{we}, the terms in the first part of the right hand side
can be estimated by
\begin{eqnarray*}
&&C\langle t \rangle^{-3/2}
\|\langle  c_I t-r\rangle^{1/2} \partial \Ga^\be u^I
\langle c_J t-r\rangle\partial^2 \Ga^\ga u^J
\partial\Gamma^\al u^K\|_{L^1(r\le c_0 t/2 )}\\
&\le& C \langle t \rangle^{-3/2}
\|\langle  c_I t-r\rangle^{1/2}\partial \Ga^\be u^I\|_{L^\infty}
\|\langle c_J t-r\rangle\partial^2 \Ga^\ga u\|_{L^2}
\|\partial \Gamma^\al u\|_{L^2}\\
&\le& C \langle t \rangle^{-3/2}
\Big[E_{|\be|+3}^{1/2}(u(t))+\mathcal{X}_{|\be|+3}(u(t))\Big]
\mathcal{X}_{|\ga|+2}(u(t))E_\eta^{1/2}(u(t))\\
&\le& C\<t\>^{-3/2} E_{\eta}(u(t))E_\od^{1/2}(u(t)) .
\end{eqnarray*}
For the terms in the second part, using \eqref{ellinf8} instead of
\eqref{ellinf6} and assuming $|\be|\le |\ga|$ (and so $|\be|\le [\eta-1]/2\le \eta-3$), we see
\begin{eqnarray*}
&&\|\partial\Gamma^\be u^I \partial\Gamma^\ga u^J
\partial\Gamma^\al u^K\|_{L^1(r\le c_0 t/2)}\\
&\le& \|\partial\Gamma^\be u^I \partial\Gamma^\ga u^J\|_{L^2(r\le c_0 t/2)}
\|\partial\Gamma^\al u^K\|_{L^2}\\
&\le& \<t\>^{-3/2}\|r\<c_I t-r\>^{1/2} \partial\Gamma^\be u^I \|_{L^\infty} \left\|
\frac{\<c_J t-r\>}{r}
\partial\Gamma^\ga u^J\right\|_{L^2} E_\eta^{1/2}(u(t))\\
&\le&\<t\>^{-3/2}\left(E^{1/2}_{|\be|+2}(u(t))+\X_{|\be|+3}(u(t))\right)\left(E^{1/2}_{|\ga|+1}(u(t))+\X_{|\ga|+2}(u(t))\right)E_\eta^{1/2}(u(t))\\
&\le&\<t\>^{-3/2}E_\eta(u(t))E_\od^{1/2}(u(t))
,
\end{eqnarray*}
where in the third inequality we have used the Hardy inequality.

The estimates for $II_j(r\le c_0 t/2)$ proceeds similarly with the obvious modification by using the Hardy inequality for terms involving $\Ga^\al u$ (and $|h_k(x)|\le Cr^{-1/2}\<r\>^{-1/2}\le C/r$).
This gives us the required upper bound for the portion of the integrals over $r\le  c_0 t/2$ in
\eqref{en3}.

{\bf Away from the origin}.
It remains to give the estimate for $ r \ge c_0 t/2$. It is here, finally, where
the difference of speed $c_I$ and the null condition enters.

{\bf Non-resonance}.
Let us start with non-resonant terms, that is, those for which
 $(I,J,K)\ne(K,K,K)$.  In this case, for $I_j+II_j(r\ge c_0 t/2)$, we need to control the quasilinear terms
\beeq\label{NR-QLterm} \sum_{(I,J)\ne(K,K)}\sum_{\be+\ga\le\al, |\al|\le \eta-1,  \ga\le \eta-2}
\left\|\partial \Ga^\be u^I \partial^2 \Ga^\ga u^J \left(|\partial\Gamma^\al u^K|+\frac{|\Gamma^\al u^K|}{r^{1/2}\<r\>^{1/2}} \right)\right\|_{L^1(r\ge c_0
t/2)}\eneq and the semilinear terms
\beeq\label{NR-SLterm}\sum_{(I,J)\ne(K,K)}\sum_{|\al|\le\eta-1,\be+\ga\le \al}
\left\|\partial\Gamma^\be u^I \partial\Gamma^\ga u^J
 \left(|\partial\Gamma^\al u^K|+\frac{|\Gamma^\al u^K|}{r^{1/2}\<r\>^{1/2}} \right)\right\|_{L^1(r\ge c_0 t/2)}.
\eneq
We separate two cases for \eqref{NR-QLterm}: $I\neq J$ and $I=J\neq K$. In the first case, we have $c_I\neq c_J$, and
$$\langle t \rangle^{3/2}
\le C \<r\>(\<c_I t-r\>+\<c_J t-r\>)^{1/2}\le C\< r \> \langle c_I t-r\rangle^{1/2}
\langle c_J t-r\rangle^{1/2}\ .$$
Using \eqref{ellinf6} and Hardy's inequality
we have the estimate
\begin{eqnarray*}
&&\|\partial\Gamma^\be u^I
\partial^2\Gamma^\ga u^J
\partial\Gamma^\al u^K\|_{L^1(r\ge c_0 t/2)}
+
\|\partial\Gamma^\be u^I
\partial^2\Gamma^\ga u^J
\Gamma^\al u^K/r\|_{L^1(r\ge c_0 t/2)}
\\
&\le& C\langle t \rangle^{-3/2}
\|\langle r \rangle \langle c_I t-r\rangle^{1/2}
\partial\Gamma^\be u^I\|_{L^\infty}
\|\langle c_J t-r\rangle^{1/2} \partial^2\Gamma^\ga u^J\|_{L^2}\\
&&\times\left(\|\partial\Gamma^\al u^K\|_{L^2}+\|\Gamma^\al u^K/r\|_{L^2}\right)\\
&\le& C\langle t \rangle^{-3/2}
 \Big[E_{|\be|+3}^{1/2}(u(t))+\X_{|\be|+3}(u(t))\Big]
\mathcal{X}_{|\ga|+2}(u(t))
E_{|\al|+1}^{1/2}(u(t))\\
&\le& C\langle t \rangle^{-3/2}
 E_\mu(u(t)) E_\od^{1/2}(u(t)).
\end{eqnarray*}
Otherwise, if $I=J\ne K$, using \eqref{ellinf6}  and \eqref{ellinf4}, we get
\begin{eqnarray*}
&&\|\partial\Gamma^\be u^I
\partial^2\Gamma^\ga u^I
\partial\Gamma^\al u^K\|_{L^1(r\ge c_0 t/2)}+
\|r^{-1/2}\<r\>^{-1/2}\partial\Gamma^\be u^I
\partial^2\Gamma^\ga u^I
\Gamma^\al u^K\|_{L^1(r\ge c_0 t/2)}
\\
&&\le Ct^{-1/2}\langle t \rangle^{-1}
\|\partial\Gamma^\be u^I\|_{L^2}
\|\langle c_I t-r\rangle^{1/2} \partial^2\Gamma^\ga u^I\|_{L^2}\\
&&\times\left(
\|\langle r \rangle \langle c_K t-r\rangle^{1/2}
\partial\Gamma^\al u^K\|_{L^\infty}
+\|\langle r\rangle^{1/2}
\Gamma^\al u^K\|_{L^\infty}
\right)
\\
&&\le Ct^{-1/2}\langle t \rangle^{-1} E_{|\be|+1}^{1/2}(u(t))
 \X_{|\ga|+2}(u(t))
\left[E_{|\al|+3}^{1/2}(u(t))+\X_{|\al|+3}(u(t))\right]\\
&&\le Ct^{-1/2}\langle t \rangle^{-1}
 E_\mu(u(t)) E_\od^{1/2}(u(t)).
\end{eqnarray*}

For the semilinear terms \eqref{NR-SLterm}, we have $I\ne K$ and so $c_I\ne c_K$. Assuming $c_I<c_K$, then
\begin{eqnarray*}
&&\left\|\partial\Gamma^\be u^I \partial\Gamma^\ga u^J
 \left(|\partial\Gamma^\al u^K|+\frac{|\Gamma^\al u^K|}{r^{1/2}\<r\>^{1/2}} \right)\right\|_{L^1(r\ge c_0 t/2)}\\
 &\le &
C t^{-1/2} \<t\>^{-1}\left\|\partial\Gamma^\be u^I \partial\Gamma^\ga u^J\right\|_{L^1}
 \|\<r\>^{1/2}\Gamma^\al u^K\|_{L^\infty}\\
 &&+C\<t\>^{-3/2}\left\|\partial\Gamma^\be u^I \partial\Gamma^\ga u^J\right\|_{L^1}
 \|\<r\>\<c_K t-r\>^{1/2}\partial\Gamma^\al u^K\|_{L^\infty(c_0 t/2<r<(c_I+c_K) t/2)}
\\
&&+C \<t\>^{-3/2}\left\|\partial\Gamma^\al u^K \partial\Gamma^\ga u^J\right\|_{L^1}
 \|\<r\>\<c_I t-r\>^{1/2}\partial\Gamma^\be u^I\|_{L^\infty(r\ge (c_I+c_K) t/2)}\\
&\le& C t^{-1/2}\<t\>^{-1} E_\mu(u(t)) E_\od^{1/2}(u(t))
\end{eqnarray*}
by using \eqref{ellinf6}  and \eqref{ellinf4}. The case $c_I>c_K$ can be handled the same way. This completes the estimates for non-resonant terms over $r\ge  c_0 t/2$ in
\eqref{en3}.

{\bf Resonance}. In the resonant case, we have
 $(I,J,K)=(K,K,K)$ and we will denote $u^K=u$ and $c_K=c$.
An application of Lemma \ref{ptwnc} yields the following upper bound
for $\<t\> I_j(r\ge c_0 t/2)$:
\begin{eqnarray*}
&&
\sum_{\be+\ga\le\al, |\ga|\le \eta-2, |\al|\le\eta-1}
\Big[
\|\Ga \Ga^{\be} u  \partial^2 \Ga^\ga u
\partial\Gamma^\al u\|_{L^1(r\ge c_0 t/2)}\\
&&+\|\partial\Gamma^{\be}u \partial \Ga \Ga^{\ga}u
\partial\Gamma^\al u\|_{L^1(r\ge c_0  t/2)}+\|\langle ct-r\rangle \partial\Gamma^{\be}u \partial^2 \Ga^\ga u
\partial\Gamma^\al u\|_{L^1(r\ge c_0  t/2)}\Big]\\
&&+
\sum_{\be+\ga\le\al, |\al|\le\eta-1}
\||\Ga \Ga^{\be} u  \partial \Ga^\ga u
\partial\Gamma^\al u^I|+\langle c t-r\rangle |\partial\Gamma^{\be}u \partial \Ga^\ga u
\partial\Gamma^\al u|\|_{L^1(r\ge c_0  t/2)}
.
\end{eqnarray*}
We still need to get an additional
decay factor of $t^{-1/2}$.

Since $r\ge c_0 t/2$, we have
$\langle r \rangle \ge C \langle t \rangle$.
Thus, we have using \eqref{ellinf4}
\begin{eqnarray*}
&&\|\Ga\Gamma^{\be} u  \partial^2 \Ga^\ga  u
\partial\Gamma^\al u \|_{L^1(r\ge c_0  t/2 )}\\
&\le& C \langle t \rangle^{-1/2}
\|\langle r \rangle^{1/2}\Ga\Gamma^{\be} u \|_{L^\infty(r\ge c_0  t/2 )}
\|\partial^2 \Ga^\ga  u \|_{L^2}
\|\partial\Gamma^\al u \|_{L^2}\\
&\le& C \langle t \rangle^{-1/2}
E_{|\be|+3}^{1/2}(u(t)) E_\eta(u(t))\\
&\le& C \langle t \rangle^{-1/2}
E_{\od}^{1/2}(u(t)) E_\eta(u(t)).
\end{eqnarray*}

In a similar fashion, the second term is handled using
\eqref{ellinf5}:
\begin{eqnarray*}
&&\|\partial\Gamma^{\be} u  \partial\Ga\Gamma^{\ga} u
\partial\Gamma^\al u \|_{L^1(r\ge c_0  t/2 )}\\
&\le& C \langle t \rangle^{-1}
\|\partial\Gamma^{\be} u \|_{L^2}
\|\langle r \rangle\partial\Ga \Gamma^{\ga} u \|_{L^\infty(r\ge c_0  t/2 )}
\|\partial\Gamma^\al u \|_{L^2}\\
&\le& C \langle t \rangle^{-1}
E_{|\ga|+3}^{1/2}(u(t)) E_\eta(u(t))\\
&\le& C \langle t \rangle^{-1}
E_{\od}^{1/2}(u(t)) E_\eta(u(t)).
\end{eqnarray*}

The third term is estimated using \eqref{ellinf5}
again and \eqref{we1}.
\begin{eqnarray*}
&&\|\langle ct-r\rangle \partial\Gamma^{\be} u  \partial^2 \Ga^\ga  u
\partial\Gamma^\al u \|_{L^1(r\ge c_0  t/2 )}\\
&\le& C \langle t \rangle^{-1}
\|\langle r \rangle \partial\Gamma^{\be} u \|_{L^\infty(r\ge  c_0 t/2 )}
\|\langle ct-r\rangle  \partial^2 \Ga^\ga  u \|_{L^2}
\|\partial\Gamma^\al u \|_{L ^2}\\
&\le& C \langle t \rangle^{-1}
E_{|\be|+3}^{1/2}(u(t))\mathcal{X}_{|\ga|+2}(u(t))E_\eta^{1/2}(u(t))\\
&\le& C \langle t \rangle^{-1}E_{\od}^{1/2}(u(t)) E_\eta(u(t)).
\end{eqnarray*}
For the terms arising in semilinear part, by \eqref{ellinf4}, \eqref{ellinf8} and the fact $\<r\>\ge C\<t\>$, we have
\begin{eqnarray*}
&&\|\Ga \Ga^{\be} u  \partial \Ga^\ga u
\partial\Gamma^\al u^I\|_{L^1(r\ge c_0  t/2)}
+\|\langle c t-r\rangle \partial\Gamma^{\be}u \partial \Ga^\ga u
\partial\Gamma^\al u\|_{L^1(r\ge c_0  t/2)}\\
&\le&
\<t\>^{-1/2}(\|\<r\>^{1/2}\Ga \Ga^{\be} u\|_{L^\infty}
+\|\<r\>^{1/2}\langle c t-r\rangle \partial\Gamma^{\be}u\|_{L^\infty})\| \partial \Ga^\ga u\|_{L^2}
\|\partial\Gamma^\al u\|_{L^2}\\
&\le&\<t\>^{-1/2}
(E_{|\be|+3}^{1/2}(u(t))+\X_{|\be|+3})E_\eta(u(t))\\
&\le& C \langle t \rangle^{-1/2}E_{\od}^{1/2}(u(t)) E_\eta(u(t)).
\end{eqnarray*}
To complete the proof of \eqref{en3}, we still need to give the estimate for $II_j(t)$. In the resonant situation, noting that we always use $L^2$ norm to control the terms involving $\pa\Ga^\al u$ in the proof of $I_j(t)$, it is easy to adapt the previous proof to get the required estimate for $II_j(t)$ by using Hardy's inequality and
the fact that $|h_k(x)|\le C r^{-1}$.

\subsection{Conclusion of the proof}\label{sec-conclusion}
We are now ready to complete the proof of Theorem \ref{thm}, by using \eqref{en2} and \eqref{en3}.

Recalling the definition \eqref{newEnergy1} and
\eqref{newEnergy2} for $E^\al_0$, $E_1^\al$ and $E_{2,k}^\al$, we know, by Sobolev embedding,
\begin{eqnarray*}
\sum_{|\al|=m-1}(|E_1^\al(t)|+\sup_{k\ge 0} |E_{2,k}^\al(t)|)&\le& C
\|\pa u\|_{L^\infty}E_{m}(u(t))\\
&\le& C E_3^{1/2}(u(t))E_{m}(u(t))\\
&\le&
  C_4  C_1\eps E_{m}(u(t))
\end{eqnarray*} for some $  C_4 \ge 1$, and  $1\le m\le \od$.

Based on this observation and the smallness assumption of $\eps$ and $\de$ (such that $8 C_3   C_4  C_1 (\ep+\de)\le 1$), we can easily obtain the following inequalities from
\eqref{en2} and \eqref{en3},
\begin{eqnarray}
&&LE_\od(T)+\sup_{0\le t\le T}E_\od(u(t))\label{en2-2}\\
&\le& 2C_3 E_\od(u(0))+2 C_3 \int_0^T \langle t \rangle^{-1} {E}_\eta^{1/2}(u(t))
{E}_\od(u(t)) dt\nonumber
\end{eqnarray}
and
\begin{eqnarray}
&&LE_\eta(T)+\sup_{0\le t\le T}E_\eta(u(t))\label{en3-2}\\
&\le&  2 C_3  E_\eta(u(0))+ 2 C_3  \int_0^T t^{-1/2}\langle t \rangle^{-1} {E}_\od^{1/2}(u(t))
{E}_\eta(u(t)) dt\nonumber
\end{eqnarray}
on any interval $[0,T]$ with
$$T<T_0=\sup\{T: E_\eta^{1/2}(u(t))\le 2C_1 \eps, t\in [0,T]\}\ .$$

Since $E_\eta^{1/2}(u(t))\le 2C_1\eps$, an application of the Gronwall inequality to \eqref{en2-2} gives us
\[
LE_\od(t)+{E}_\od(u(t))\le   2 C_3  e^{4 C_3 C_1\eps} {E}_\od(u(0))
 \<t\>^{4 C_3 C_1\eps}.
\]
  Inserting this bound into
\eqref{en3-2}, we obtain
\begin{eqnarray*}
&&LE_\eta(t)+{E}_\eta(u(t))\le  2 C_3  {E}_\eta(u(0))\\
&&\times\exp \left( (2 C_3)^{3/2} e^{ 2 C_3  C_1\ep}E_\od^{1/2}(u(0)) \int_0^\infty \<t\>^{-1+2 C_3 C_1\eps}t^{-1/2} dt\right)\ .
\end{eqnarray*}
Setting $C_1=\sqrt{ 2 C_3 }$, $$C_0=C_1^3 e^{C_1^3\ep}\int_0^\infty \<t\>^{-1+C_1^3\eps}t^{-1/2} dt/2\ ,$$ and $C_2=\max(C_1^2 e^{2 C_1^3\ep},2C_1^3)$, we have
$$LE_{\eta}(t)+{E}_\eta(u(t))\le  2 C_3  {E}_\eta(u(0))
\exp (2 C_0 E_\od^{1/2}(u(0)) )< 2 C_3 \eps^2=(C_1\eps)^2\ ,$$
$$LE_{\od}(t)+{E}_\od(u(t))\le C_2 {E}_\od(u(0))
\<t\>^{C_2\eps}\ .$$
To ensure the finiteness of $C_0$, in addition to the smallness assumption on $\ep$ to absorb the perturbation, we need also to require $C_1^3 \ep<1/2$.

With this we see that ${E}^{1/2}_\eta(u(t))$ remains less than
$C_1\eps$ throughout the interval $0\le t < T_0$. A standard continuity argument
shows that $E_\eta(u(t))$ is bounded for all time, which
completes the proof of Theorem \ref{thm}.

\section{Appendix: the case of asymptotically flat manifolds}\label{sec-App}
In this appendix, we will give the almost global existence (and global existence respectively) for the Cauchy problem of the quasilinear wave system  when spatial dimension is three (and higher), posed on certain asymptotically flat manifolds.

Let us begin with the space-time manifolds. We consider the asymptotically
flat Lorentzian manifolds $( \R^{1+n} , \gm)$ with
\begin{equation*}
\gm =  g_{\al\be} (t,x) \, d x^\al \, d x^\be =\sum_{\al,\be=0}^n
 g_{\al\be} (t,x) \, d x^\al \, d x^\be
.
\end{equation*}
The metric $\gm$ is assumed to be a small asymptotically flat perturbation of the Minkowski metric. More precisely, we suppose $g_{\al\be} (t,x) \in C^{\infty} ( \R^{1+n} )$ and, for some
fixed $\rho
>0$ and $\delta\ll 1$,
\begin{equation}\tag{H2} \label{a-H2}
\forall \ga \in \N^{1+n} \qquad |\partial^\ga_{t,x} ( g_{\al\be}(t,x) -
m_{\al\be} )|\le C_\ga \delta  \< x \>^{- \vert \ga \vert - \rho}  ,
\end{equation}
with $(m_{\al\be})=Diag(1,-1,-1,\cdots,-1)$ being the standard Minkowski metric and $\<x\>=\sqrt{1+|x|^2}$. An example of such metric can be $$g_{\al\be}=m_{\al\be}+\de \<x\>^{-\rho}+\de\phi(t/\<x\>) \<x\>^{-\rho}$$ with $\phi\in C_0^\infty$.
Since $\delta\ll 1$, it is clear that the metric $\gm$ is a non-trapping perturbation. Let $g =(-1)^n \det (g_{\al\be})$, the Laplace--Beltrami
operator associated with $\gm$ is given by
$$\Box_{\gm} =  \sqrt{g}^{-1} \partial_\al g^{\al\be} \sqrt{g} \partial_\be ,
$$
where $(g^{\al\be} (t,x))$ denotes the inverse matrix of $(g_{\al\be}(t,x))$.

We would also like to investigate the case of radial metric, by which we mean that, when writing out the metric in polar coordinates $(t,x)=(t,r\omega)$ with $\omega\in\mathbb{S}^{n-1}$, we have
$$\gm=\tilde g_{00}(t,r)dt^2+2\tilde g_{01}(t,r)dtdr+\tilde g_{11}(t,r)dr^2+\tilde g_{22}(t,r) r^2 d\omega^2\ .$$ In this form, the assumption \eqref{a-H2} on asymptotic flatness is equivalent to the following requirement
\begin{equation}\tag{H2'} \label{a-H2'}
|\partial^\ga_{t,x} (\tilde g_{00}-1, \tilde g_{11}+1,\tilde g_{22}+1,\tilde g_{01})|
\le C_\ga \delta  \< x \>^{- \vert \ga \vert - \rho}  .
\end{equation}

Consider the initial value problem for the quasilinear wave equations of the form
\begin{equation}
\label{a-pde}
(\Box_\gm u)^I = N^I(u,u), I=1,2,\cdots, M
\end{equation}
in which the  quadratic nonlinearity $N=Q+S$ is of the form \eqref{non}.
Our construction of solutions will depend on the energy integral method, which requires the quasilinear part to be symmetric
\eqref{symmetric}.

In contrast to the null-form system, as in \cite{MeSo06}, we will be able to avoid the use of the scaling vector field $S$, and the vector fields to be used will be labeled as
\[
Y=(Y_0,\ldots,Y_6)=(\partial,\Omega).
\]

For the energy norm, we will use the standard energy norm
$$E_1(u(t))=\frac{1}2\sum_{I=1}^M\int_{\R^n} |\pa u^I(t,x)|^2 dx\ .
$$
The higher order derivatives will be estimated through
\begin{equation}
\label{a-ennorm}
E_m(u(t))=\sum_{|\al|\le m-1} E_1(Y^\al u(t)),
\qquad m=2, 3, \cdots\ .
\end{equation}

As for the null-form systems, an important intermediate role will be played by the following local energy norm
$$\I_m(u(t))= \sum_{I=1}^M\sum_{|\al|\le m-1}
\left\|r^{-1/2+\mu}\<r\>^{-\mu'} \left(|\partial \Gamma^\al u^I (t)|+\frac{|\Gamma^\al u^I (t)|}r\right)\right\|_{L^2(\R^n)}^2$$
with $\mu\in (0,1/2)$ and $\mu'> \mu$ to be determined later.
This norm is extracted from the local energy norm, which is defined as
\begin{equation}
\label{a-LEnorm}
LE_m(t)= \int_0^t \I_m(u(\tau))d\tau\ .
\end{equation}
 We will choose $\mu=1/4$.
When $n\ge 4$, the choice for $\mu'$ will be $\min(n-2,2\rho-1,3)/4$ for the general case and $\mu'=\min(n-2,1+2\rho,3)/4$
for the radial metric. In the case of $n=3$, we will set $\mu'=\min(2\rho-1,3)/4$ (and $\mu'=\min(2\rho+1,3)/4$ for the radial metric).

In order to describe the solution space, we introduce
\[
H^m_Y(\R^n)=\{f\in L^2(\R^n) : (\nabla,\Omega)^\al f\in L^2,
\; |\al|\le m\},
\]
with the norm
\begin{equation}
\label{a-sobnorm}
\|f\|_{H^m_Y}=\sum_{|\al|\le m}\| (\nabla,\Omega)^\al f\|_{L^2}.
\end{equation}
Solutions will be constructed in the space
$\dot{H}_Y^m(T)$ obtained by closing
the set $C^\infty([0,T);C_0^\infty(\R^n))$
in the norm $\sup\limits_{0\le t< T} E_m^{1/2}(u(t))$.
Thus,
$$
\dot{H}_Y^ m(T)\subset\left\{u(t,x) : \partial u(t,\cdot)
\in\bigcap_{j=0}^{ m-1} C^j([0,T);H^{ m-1-j}_Y)\right\}.
$$
By Sobolev embedding, it follows that
$\dot{H}_Y^m(T)\subset C^{m-[(n+2)/2]}([0,T)\times\R^n)$.

Let us now state our main result precisely.
\begin{theorem}
\label{a-thm}
Let $n\ge 3$, $\de\ll 1$, $\rho>1$ for the general metric and $\rho>0$ for the radial metric.
Assume that the nonlinear terms in \eqref{a-pde}
satisfy the symmetric condition \eqref{symmetric}.
Then there exist constants $\eps_0, c_0\ll 1$, such that the Cauchy problem for \eqref{a-pde} has a unique global (almost global for $n=3$) solution $u\in \dot H^\od_Y(t)$ for $t\in [0,T_\ep]$  with
$$T_\ep=\left\{
\begin{array}{ll}
\infty&n\ge 4\\
\exp(c_0/\ep)& n=3
\end{array}
\right.,$$ when the initial data
satisfy
\begin{equation}
\label{a-smallness}
E_{\od}^{1/2}(u(0))= \ep\le \ep_0, \od\ge n+4.
\end{equation}
Moreover, the solution satisfies the bounds for some $C_1>1$,
\[
\sup_{t\in[0,T_\ep]}E_{\od}(u(t))+LE_\od (T_\ep)+\de_{3n}\frac{\ep}{2c_0} KSS_\od(T_\ep)\le C_1^2 \ep^2.
\]
\end{theorem}

\subsection{Commutation with vector fields}
In preparation for the energy estimates, we need to consider
the commutation properties of the vector fields $Y$ with respect
to the nonlinear terms.

\begin{lemma}
\label{a-nd}
Let $u$ be a solution of \eqref{a-pde}.
Assume that the nonlinearity is of the form
 \eqref{non}.
Then for any $\al\in\N^7$,
$$
\Box_\gm Y^\al u = \sum_{\be+\ga+\mu=\al}N_\mu^\al(Y^\be u,Y^\ga u)+\sum_{|\be|\le|\al|-1}\(r_0 \nabla^2 \Gamma^\be u+r_1 \nabla \Gamma^\be u\right),
$$
in which each $N_\mu^\al$ is a quadratic nonlinearity of the form
\eqref{non}, and $r_m$ with $m\in\N$ denote functions such that
$$|\pa^\al r_m(t,x)|\le C_\al \delta \<r\>^{-\rho-m-|\al|}\quad\mbox{for any}\quad\al\in\N^{1+n}\ .$$
Moreover, if $|\mu|=0$, then $N_\mu^\al=N$.
In addition, if the metric is radial,
$$
\Box_\gm Y^\al u = \sum_{\be+\ga+\mu=\al}N_\mu^\al(Y^\be u,Y^\ga u)+\sum_{|\be|\le|\al|-1}\(r_1 \nabla^2 \Gamma^\be u+r_2 \nabla \Gamma^\be u\right) .
$$
\end{lemma}

\begin{prf}
It is easy to check that
\[
[Y,N](u,v)=Y N(u,v) - N(Y u,v)-N(u,Y v)
\]
 is a quadratic nonlinearity of the form \eqref{non}.

By \eqref{a-H2}, we have
$$\Box_\gm =\Box+r_0\pa^2+r_1\pa\ ,$$
with $\Box=Diag(\pt^2- \Delta,\cdots, \pt^2-\Delta)$.
Recall that
$$[\Box, Y_j]=0\ .$$
We want to prove the result by induction. It is clear the result is true for $|\al|=0$. Now assume that it is true for any $\al$ with $|\al|=m$. Given $\al_0$ with $|\al_0|=m+1$, we can find some $j$ and $\al$ with $|\al|=m$ and $Y^{\al_0}=Y_jY^\al$. Then  by the inductional assumption, we can calculate $ \Box_\gm Y^{\al_0} u=\Box_\gm Y_j Y^\al u$ as follows
\begin{eqnarray*}
   &&[\Box_\gm,Y_j] Y^\al u+Y_j \Box_\gm Y^\al u\\
  &=&[r_0\pa^2+r_1\pa,
  Y_j]Y^\al u\\
  &&+ \sum_{\be+\ga+\mu=\al}Y_j N_\mu^\al(Y^\be u,Y^\ga u)+\sum_{|\be|\le|\al|-1}Y_j\(r_0 \pa^2 Y^\be u+r_1 \pa Y^\be u\right)\\
  &=&\sum_{|\be|\le|\al|}\(r_0 \pa^2 Y^\be u+r_1 \pa Y^\be u\right)\\
  &&+ \sum_{\be+\ga+\mu=\al}\left\{[Y_j, N_\mu^\al](Y^\be u,Y^\ga u)+N_\mu^\al(Y_jY^\be u,Y^\ga u)+N_\mu^\al(Y^\be u, Y_jY^\ga u)\right\}
\end{eqnarray*}
which is of the required form.  This completes the proof for the general case. When the metric is radial,
the rotational vector fields $\Omega$ are commutative with $\Box_\gm$ and we need only to give the estimate for $\pa$. In this case, we have $[r_{2-j} \pa^j,\pa]=r_{3-j} \pa^j$ with $j=1,2$ and the same argument will give the proof.
\end{prf}
\subsection{Sobolev-type inequalities}
\label{a-sob}
As in \cite{MeSo06}, we need to use the Sobolev inequalities which do not involve the Lorentz boost operators and scaling vector field, which is also related with the trace estimates.

Recall that we have the following trace estimate (see e.g. (1.3) in \cite{FaWa})
$$r^{s}\|f(r\omega)\|_{H_\omega^{(n-1)/2-s}}\le C \|f\|_{\dot H^{n/2-s}}, 0<s<(n-1)/2$$
which gives us  (with $s=1/4$)
\beeq\label{a-120703-1}r^{1/4}\|f(r\omega)\|_{L^\infty_\omega}\le C \sum_{|\al|+|\be|\le (n+3)/2}\|\pa^\al\Omega^\be f\|_{L^2}.\eneq
Moreover, Sobolev embedding in the polar coordinates gives us for $r\ge 2$ (see e.g. (2.1) in \cite{MeSo06})
\beeq\label{a-120703-2}r^{(n-1)/2}\|f(r\omega)\|_{L_\omega^{\infty}}\le C \sum_{|\al|+|\be|\le n/2+1}\|\pa^\al\Omega^\be f\|_{L^2(r-1<|x|<r+1)}\ .\eneq
As a direct consequence of the estimates \eqref{a-120703-1} and \eqref{a-120703-2}, we have
\begin{eqnarray}\label{a-120703-3}
&&\| f g \|^2_{L^2}
\\
&\le &\|fg\|_{L^2(|x|\le 2)}^2+\sum_{j\ge 3}\|fg\|_{L^2(j-1\le |x|\le j+1)}^2\nonumber\\
&\le &\|r^{1/4}f\|^2_{L^\infty(|x|\le 2)}\|r^{-1/4}g\|^2_{L^2(|x|\le 2)}\nonumber\\
&&+\sum_{j\ge 3}
\|r^{(n-1)/2}f\|_{L^\infty(j-1\le |x|\le j+1)}^2\|r^{-(n-1)/2} g\|_{L^2(j-1\le |x|\le j+1)}^2\nonumber
\\
&\le& C\sum_{|\al|\le (n+3)/2}\|Y^\al f\|_{L^2(|x|\le 3)}^2\|r^{-1/4}g\|_{L^2(|x|\le 2)}^2\nonumber\\
&&+C\sum_{j\ge 3}
\sum_{|\al|\le (n+3)/2}\|Y^\al f\|_{L^2(j-2\le |x|\le j+2)}^2
\|r^{-(n-1)/2} g\|_{L^2(j-1\le |x|\le j+1)}^2\nonumber
\\
&\le&
C\|r^{-1/4}\<r\>^{-(n-2)/4}g\|_{L^2}^2\sum_{|\al|\le (n+3)/2}\|r^{-1/4}\<r\>^{-(n-2)/4}Y^\al f\|_{L^2}^2
\nonumber
\end{eqnarray}

\subsection{Energy and local energy estimates}\label{a-local}
As in the beginning of Section \ref{sec-energy},
we assume that $u(t)\in \dot{H}^\od_Y(T)$ with $k\ge n+4$ is a local
solution of the initial value problem for \eqref{a-pde}.
Our task will be to show that $E_\od(u(t))$ remains finite
for all $t\ge0$ when $n\ge 4$ (and for some interval $[0,\exp(c/\ep)]$ with $c\ll 1$ when $n=3$).  To do so, we will derive integral inequalities for $E_\od(u(t))+LE_\od(t)$.
If \eqref{a-smallness} holds, then $E_\od^{1/2}(u(0))=\ep$ and $E_\od(u(t))
+LE_\od(T)+(\log(2+T))^{-1}KSS_\od(T)
< (C_1\eps)^2$ for certain $C_1>1$ and small interval $t\in [0,T]$.
Define
$$T_0=\sup\{T: E_\od(u(t))
+LE_\od(T)+(\log(2+T))^{-1}KSS_\od(T)
\le (2 C_1\eps)^2, t\in [0,T]\}\ .$$
Here the constant $C_1$ will be determined later.
All of the following computations will be valid on
the time interval $[0,T_0)$.

It is easy to check that the same argument of Section \ref{sec-energy} gives us the required local energy estimates, that is, by writting
$$\Box_\gm= \Box+ (g^{\al\be}-m^{\al\be})\pa_\al\pa_\be+r_1 \pa \ ,$$
and applying Lemma \ref{KSS} with $h^{\al\be}=g^{\al\be}-m^{\al\be}$, Lemma \ref{a-nd} and the symmetry condition \eqref{symmetric},  we have
\begin{eqnarray}
&&\sup_{t\in[0,T]}E_m(u(t))+LE_m(T)+(\log(2+T))^{-1}KSS_m(T)\nonumber\\
&\le& C E_m(u(0))+C \de LE_m(T)\nonumber\\
&&+C\sum_{|\al|\le m-1} \left|\int_0^T C_1^{\al}(t)dt\right| +C \sum_{|\al|\le m-1} \sup_{k\ge 0}\left|\int_0^T C_{2,k}^{\al}(t)dt\right|\nonumber\\
&\le & C E_m(u(0))+C C_1(\de+\ep) LE_m(T)\nonumber\\
&&+ C\sum_{|\al|=m-1} |E_1^{\al}(T)-E_1^{\al}(0)|
+ C\sup_{k\ge 0}\sum_{|\al|=m-1} |E_{2,k}^{\al}(T)-E_{2,k}^{\al}(0)|
\nonumber\\
&&+C\sum_{I=1}^M \sum_{|\al|\le m-1, \be+\ga+\mu=\al} \sup_{k\ge 0}\left|\int_0^T\int_{\R^3} S_\mu^{\al,I}(\Gamma^\be u,\Gamma^\ga u) (\pt,L_k) \Ga^\al u^I(t) dx dt\right|
\nonumber\\
&&+C \sum_{I=1}^M \sum_{|\al|\le m-1, \be+\ga+\mu=\al, |\ga|<m-1} \sup_{k\ge 0}\left|\int_0^T\int_{\R^3} Q_\mu^{\al,I}(\Gamma^\be u,\Gamma^\ga u) (\pt,L_k) \Ga^\al u^I(t) dx dt\right|
\nonumber\\
&&+C \sum_{I=1}^M \sum_{|\al|=m-1}\sup_{k\ge 0}\left|\int_0^T \int_{\R^3}Q^{I,\be\mu\ga}_{JK}   \pa_\ga \pa_\be u^J \pa_{\mu} \Ga^\al u^K (\pt,L_k) \Ga^\al u^I dxdt \right|\nonumber\\
&&+  C\sum_{I=1}^M  \sum_{|\al|=m-1}\sup_{k\ge 0}\left|\int_0^T \int_{\R^3} Q^{I,\be\mu\ga}_{JK} (\pt,h^i_k\pa_i)\pa_\be u^J \pa_{\mu} \Ga^\al u^K  \pa_\ga  \Ga^\al u^I dx dt \right|\ ,\nonumber\end{eqnarray}
where we have used the assumption
\beeq\label{a-120630-1}2\mu'\le \rho-1+2\mu,\ \rho>1\ \eneq
for general metric.
When the metric is radial, the assumption can be weakened to be
\beeq\label{a-120630-1-rad}2\mu'\le \rho+2\mu,\ \rho>0 \ (\gm \textrm{ radial})\ .\eneq
Here, as before,
$L_k=f_k 	(\pa_r +\frac{1}{r} )=h_k^i(x)\pa_i+ h_k(x)$ with
$h_k(x)=f_k/r$ and $h_k^i(x)=f_k x^i/r$, $\mu=1/4<\mu'\le 3/4$,
\beeq\label{a-newEnergy1}E_1^{\al}(t)=
\sum_{I=1}^M (Q^{I,\be\mu 0}_{JK}\delta_0^\ga-\frac 12 Q^{I,\be\mu\ga}_{JK}) \int_{\R^3} \pa_\be u^J \pa_{\mu} \Ga^\al u^K
\pa_\ga\Ga^\al u^I  dx,\eneq
and \beeq\label{a-newEnergy2}
 E_{2,k}^{\al}=\sum_{I=1}^M
Q^{I,\be\mu 0}_{JK} \int_{\R^3} \pa_\be u^J \pa_{\mu} \Ga^\al u^K  L_k \Ga^\al u^I dx
. \eneq

\subsection{Conclusion of the proof}
We are now ready to complete the proof of Theorem \ref{a-thm}, by using \eqref{a-120703-3} and the local energy estimates in the Subsection \ref{a-local}.

We need only to give a better upper bound for
$$\sup_{t\in[0,T]}E_\od(u(t))+LE_\od(T)+(\log(2+T))^{-1}KSS_\od(T)\ .$$
By the local energy estimates in the Subsection \ref{a-local}, we see that there exists a universal constant $C_2\ge 1$ such that
\begin{eqnarray*}
&&\sup_{t\in[0,T]}E_\od(u(t))+LE_\od(T)+(\log(2+T))^{-1}KSS_\od(T)\label{0610-1}\\
&\le& C_2 E_\od(u(0))+C_2C_1(\de+\ep) LE_\od(T)+C_2\sum_{|\al|=\od-1} \sup_{k\ge 0}|E_{2,k}^{\al}(T)-E_{2,k}^{\al}(0)|\nonumber\\
&&+C_2\sum_{|\al|=\od-1} |E_{1}^{\al}(T)-E_{1}^{\al}(0)|
\nonumber\\
&&+C_2\sum_{|\al|\le\od-1}\left(\|\pa Y^\al u\|_{L^\infty_t L^2_x ([0,T]\times\R^n)}+
\left\|\frac{1}r Y^\al u\right\|_{L^\infty_t L^2_x ([0,T]\times\R^n)}\right)\nonumber\\
&&\times\sum_{|\al|\le\od-1,|\be|\le (\od-1)/2}\|\pa Y^\be u \pa Y^\al u\|_{L^1_t L^2_x ([0,T]\times\R^n)}\nonumber\\
&\le& C_2 E_\od(u(0))+C_2C_1(\de+\ep) LE_\od(T)+C_2\sum_{|\al|=\od-1} \sup_{k\ge 0}|E_{2,k}^{\al}(T)-E_{2,k}^{\al}(0)|\nonumber\\
&&+C_2\sum_{|\al|=\od-1} |E_{1}^{\al}(T)-E_{1}^{\al}(0)|
\nonumber\\
&&+C_2\sup_{t\in [0,T]}E_\od^{1/2}(u(t))\sum_{|\al|\le\od-1, (\od+n+2)/2}\|r^{-1/4}\<r\>^{-(n-2)/4}\pa Y^\al u\|_{L^2_t L^2_x ([0,T]\times\R^n)}^2,\nonumber
\end{eqnarray*}
where in the last inequality, we have used the Hardy inequality and \eqref{a-120703-3}.

Recalling the definition \eqref{a-newEnergy1} and
\eqref{a-newEnergy2} for $E_1^\al$ and $E_{2,k}^\al$, we know, by Sobolev embedding,
\begin{eqnarray*}
\sum_{|\al|=m-1}(|E_1^\al(t)|+\sup_{k\ge 0} |E_{2,k}^\al(t)|)&\le& C
\|\pa u\|_{L^\infty}E_{m}(u(t))\\
&\le& C E_3^{1/2}(u(t))E_{m}(u(t))\\
&\le&
  C_3  C_1\eps E_{m}(u(t))
\end{eqnarray*} for some $  C_3 \ge 2$.

Based on this observation and the smallness assumption of $\eps$ and $\de$ (such that $8 C_2   C_3C_1 (\ep+\de)\le 1$), recalling also that $(\od+n+2)/2\le \od-1$ since $\od\ge n+4$,
we can easily obtain the following inequality, 
\begin{eqnarray}
&&\sup_{t\in[0,T]}E_\od(u(t))+LE_\od(T)+(\log(2+T))^{-1}KSS_\od(T)\label{a-en2-2}\\
&\le& 2C_2 E_\od(u(0))+2 C_2 E_\od^{1/2}(\log(2+T))^{\de_{3n}}(LE_\od(T)+(\log(2+T))^{-1}KSS_\od(T))\nonumber
\end{eqnarray}
on any interval $[0,T]$ with
$T<T_0$.

If $n\ge 4$, we do not have the term $\log(2+T)$ and
\beeq\label{120703-f}\sup_{t\in[0,T]}E_\od(u(t))+LE_\od(T)+(\log(2+T))^{-1}KSS_\od(T)\le 3 C_1^2\ep^2\eneq
by setting
 $C_1=\sqrt{C_2}$ and $\de, \ep \le \ep_0=1/(16C_3 C_2C_1)$.
 In the case of $n=3$, we have \eqref{120703-f} for $T\le T_\ep$ and $T<T_0$, if we set
 $C_1=\sqrt{C_2}$,
 $\de, \ep\le\ep_0=1/(16C_3C_2C_1)$, and
 $$T_\ep=\exp\left(\frac{1}{16 C_2 C_1 \ep}\right)-2 \ge \exp(c_0/\ep) $$
 with $c_0=1/(32 C_2 C_1)$. When $n\ge 4$, we denote $T_\ep=\infty$. Here, we may need to let $\de$  be even smaller such that the assumption in Lemma \ref{KSS} is satisfied.

With this we see that the sum of ${E}^{1/2}_\od(u(t))$ and local energy norms remain less than
$\sqrt{3}C_1\eps$ throughout the interval $0\le t < T_0$ (and $t< T_\ep$). A standard continuity argument
shows that $E_\od(u(t))$ is bounded for $t\in [0,T_\ep)$, which
completes the proof of Theorem \ref{a-thm}.

\medskip 
{\bf Acknowledgements.}
The authors would like to thank Jason Metcalfe and Chris Sogge for helpful discussions on problems in connection with this work.


\end{document}